\documentclass[amstex]{amsart}
\usepackage{amsmath}
\usepackage{amsthm}
\usepackage{amssymb}
\usepackage{amsfonts}
\usepackage{amscd}
\theoremstyle{plain}
	\newtheorem{thm}{Theorem}[section]
	\newtheorem*{thm*}{Theorem}
	\newtheorem{cor}{Corollary}[section]
	
	\newtheorem{prop}[thm]{Proposition}
	
	\newtheorem{conj}{Conjecture}[section]
\theoremstyle{definition}
	\newtheorem{defn}{Definition}[section]
	\newtheorem{exmp}{Example}[section]
\theoremstyle{remark}
	\newtheorem*{rem}{Remark}

\newcommand{\Z}{{\bf Z}}

\newcommand{\Q}{{\bf Q}}
\newcommand{\C}{{\bf C}}
\newcommand{\BP}{{\bf P}}

\newcommand{\cO}{{\mathcal O}}

\newcommand{\cF}{{\mathcal F}}

\newcommand{\cG}{{\mathcal G}}

\newcommand{\M}{\overline{M}}

\newcommand\Spec{\mathop{\rm Spec}\nolimits}
\newcommand\cHom{\mathop{{\mathcal H}om}\nolimits}
\newcommand\cExt{\mathop{{\mathcal E}xt}\nolimits}
\newcommand\Ext{\mbox{\rm Ext}}

\newcommand\Pic{\mathop{\rm Pic}\nolimits}

\newcommand\lra{\longrightarrow}

\def\C{{\mathbb C}}
\def\F{{\mathcal F}}
\def\H{{\mathbb H}}

\def\P{{\mathbb P}}
\def\BP{{\mathbb P}}
\def\Q{{\mathbb Q}}

\def\Z{{\mathbb Z}}
\def\E{{\mathcal E}}
\def\M{{\mathcal M}}
\def\O{{\mathcal O}}
\begin{document}
\title{Relative Lefschetz Action and  BPS State Counting}
\author{Shinobu Hosono}
\address{Graduate School of Mathematical Sciences,  University of Tokyo, Komaba 3-8-1, Meguro-ku, Tokyo 153, Japan  }
\email{hosono@ms.u-tokyo.ac.jp} 
\author{Masa-Hiko Saito}
\address{Department of Mathematics, Faculty of Science,
Kobe University, Rokko, 657-8501, Kobe, Japan}
\email{mhsaito@math.kobe-u.ac.jp} 
\thanks{Partly supported by Grant-in Aid
for Scientific Research (B-09440015), (B-12440008) and (C-11874008), 
the Ministry of
Education, Science and Culture, Japan }
	\author{Atsushi Takahashi}
\address{Research Institute for Mathematical Sciences, Kyoto University, Kyoto 606-8502, Japan}
\email{atsushi@kurims.kyoto-u.ac.jp}

\begin{abstract}
In this paper, we propose a mathematical definition of  a  new 
``numerical invariants" of Calabi--Yau 3-folds from stable sheaves of dimension one, which is motivated by the Gopakumar-Vafa conjecture \cite{gv:1} in M-theory.  
Moreover, we show that  for any projective morphism 
$f:X \lra Y$ of normal projective varieties, there exists  
a natural $sl_2 \times sl_2$ action on the intersection 
cohomology group $IH(X, \Q)$ which fits into the 
perverse Leray spectral sequence. 

\end{abstract}

\maketitle
%
%
%%%%%%%%%%%%%%%%%%%%%%%%%%%%%%%%%%%%%%%%%%%%%%%%%%%%%%%%%%%%%%%%%%%%%%%%%%%%%%
%%%%%%%%%%%%%%%%%%%%%%%%%%%%%%%%%%%%%%%%%%%%%%%%%%%%%%%%%%%%%%%%%%%%%%%%%%%%%%
\section{Introduction}
~~~

Let $X$ be a Calabi--Yau $3$-fold with $\pi_1(X)=\{1\}$ and 
let us fix an ample line bundle $\O_X(1)$ on $X$. 

For $ \beta \in H_2(X, \Z)$ and an  integer $g \geq 0$, 
we denote the 0-point genus $g$ Gromov--Witten invariants of $X$ 
in the homology class $\beta$ by 
$$
N_g(\beta):=[\overline{\M_{g,0}}(X,\beta)]^{virt}\in 
A_0(\overline{\M_{g,0}}(X,\beta))\simeq\Q.
$$
Recently, physicists Gopakumar and Vafa \cite{gv:1} introduced the following remarkable 
formula for the generating function of Gromov--Witten invariants 
based on the string duality between Type IIA and M-theory.
\begin{conj}\label{conj:gvc} $($\cite{gv:1}$)$ \ \
There should exist {\bf integers} $n_h(\beta)$ such that
\begin{equation}\label{eq:gv}
\boxed{ \ \ 
\sum_{g\ge 0, \beta\in H_2(X,\Z)}N_g(\beta)q^\beta\lambda^{2g-2}=
\sum_{k>0,h\ge 0, \beta\in H_2(X,\Z)}n_h(\beta)\frac{1}{k}
\left(2\sin(\frac{k\lambda}{2})\right)^{2h-2}q^{k\beta},  \ \ }
\end{equation}
where $q^\beta:=\exp(-2\pi\int_\beta\omega)$, $\omega:=c_1(\O_X(1))$.
\end{conj}
They also proposed that integers $n_h(\beta)$ should be defined by the spin contents of 
the BPS states of M2-branes wrapped around the curves in $X$.
More precisely, they expect that a suitable D-brane moduli space $M_\beta$ and 
the natural support map $\pi_\beta:M_\beta\to S_\beta$ exist. 
By assuming the existence of an $(sl_2)_L\times (sl_2)_R$-action on some suitable 
cohomology group $H^*(M_\beta)$, 
they decompose $H^*(M_\beta)$ and rearrange it as
$$
H^*(M_\beta)=\bigoplus_{h\ge 0}\left[(\frac{1}{2})_L\oplus 2(0)_L\right]^{\otimes h}
\otimes R_h(\beta),
$$
to define numerical invariants
$$
n_h(\beta):=Tr_{R_h(\beta)}(-1)^{2H_R}.
$$
To complete their conjecture we need to define mathematically their integral 
``numerical invariants" of Calabi--Yau 3-folds by the moduli space of 
``D-branes" and to formulate their conjecture as an equivalence of the new invariants 
and Gromov--Witten invariants.
For this purpose, we have to
\begin{enumerate}
\item define the moduli space of D-branes,
\item prove the existence of an $(sl_2)_L\times (sl_2)_R$-action on a suitable cohomology 
on the above moduli space,
\item prove the Gopakumar--Vafa formula \eqref{eq:gv}.
\end{enumerate}
In this paper we deal with the first two steps and present nontrivial 
evidences for Gopakumar--Vafa conjecture. 
Especially, we can provide the answer of the problem (ii) using the intersection 
cohomology of the D-brane moduli spaces and the decomposition theorem due to \cite{bbd:1}. 
As for the D-brane moduli space we propose a natural definition in section 3.
Here is the brief plan of this paper.
In Section 2, we recall the general theory of perverse sheaves,
and prove that for any projective morphism $f:X \to Y$ there exists a natural 
$(sl_2)_L\times (sl_2)_R$-action on intersection cohomology $IH^*(X)$.
In Section 3 we consider a suitable moduli space $M_\beta$ of semi-stable sheaves on 
Calabi--Yau $3$-folds and define the support map $\pi_\beta:M_\beta\to S_\beta$.
We propose $M_\beta$ as the moduli space of D-branes and 
by applying the results in Section 2 to $\pi_\beta$, 
we obtain the numerical invariants $n_h(\beta)$.
In Section 4, we provide some evidences for Gopakumar--Vafa conjecture.

In \cite{bp:1}, Bryan and Pandharipande  proved the integrality of BPS invariants $n_h(\beta)$ coming from the formula \eqref{eq:gv} for 
Gromov--Witten invariants of some super-rigid curves in 
a Calabi--Yau 3-fold   by evaluating the virtual fundamental classes.
Also, we are informed that Fukaya--Ono \cite{fo}  proved the genus $0$ part of this conjecture in the symplectic category. 
%%%%%%%%%%%%%%%%%%%%%%%%%%%%%%%%%%%%%%%%%%%%%%%%%%%%%%%%%%%%%%%%%%%%%%%%%%%%%%
%%%%%%%%%%%%%%%%%%%%%%%%%%%%%%%%%%%%%%%%%%%%%%%%%%%%%%%%%%%%%%%%%%%%%%%%%%%%%%

\section{Relative Lefschetz Action, --- General Theory by BBD}

In this section, we  recall 
 the definition of perverse sheaves and the intersection cohomology briefly. 
Since we only need the formulation of perverse sheaves and the relative hard 
Lefschetz theorem for perverse direct image sheaves (\ref{thm:relative}), 
we shall not give any proof. 
For detail, we refer the readers to  \cite{bbd:1}.

Let $X$ be a normal complex algebraic variety.
As in [2.2.1, \cite{bbd:1}], we will only 
consider  stratifications $X = \coprod_{i=1}^rX_i$ by equidimensional algebraic strata $X_i$.

\begin{defn}$(${\bf Constructible Sheaves}$)$([2.2.1, \cite{bbd:1}]). \\ 
A $\C_X$-module ${\mathcal F}$ is called {\it constructible} if 
there exists a stratification $X=\coprod_{i=1}^rX_i$ such that 
restrictions ${\mathcal F}|_{X_i}$ are local systems on $X_i$.
\end{defn}

Let $D^b_c(\C_X)$ be the derived category of bounded complexes of $\C$-modules 
with constructible cohomology sheaves.
Let $X = \coprod_{i=1}^r X_i$ be a stratification of $X$.  In order to define perverse sheaves, 
we have to fix a perversity $p$.  
As in [2.1.16, \cite{bbd:1}], it is convenient to take the auto-dual perversity which is defined 
for each strata $ j:S \hookrightarrow X $ as 
$$
p(S) = - \dim_{\C} S. 
$$

\begin{defn}$(${\bf Perverse sheaves}$)$\\
A {\it perverse} $\C_X$-module (with the middle perversity) is an object $K^\bullet\in D^b_c(\C_X)$ such 
that the following conditions are satisfied:
\begin{enumerate}
\item (Support condition)
$$
\dim_\C{\rm supp}H^i(K^\bullet)\le -i,~~~i\in\Z.
$$
\item (Support condition for Verdier dual)
$$
\dim_\C{\rm supp}H^i({\mathbb D}_X K^\bullet)\le -i,~~~i\in\Z,
$$
where ${\mathbb D}_X$ is a Verdier dualizing functor.
Let ${}^pD^{\le 0}(\C_X)$ (resp. ${}^pD^{\ge 0}(\C_X)$) be 
the subcategory of $D^b_c(\C_X)$ whose objects are complexes 
$K^\bullet\in D^b_c(\C_X)$ satisfying the support condition 
(resp. support condition for Verdier dual).  
Let us set 
$$
Perv(\C_X):={}^pD^{\le 0}(\C_X)\cap{}^pD^{\ge 0}(\C_X).
$$
\end{enumerate}
\end{defn}
\begin{rem}
The category of perverse $\C_X$-modules is an abelian category 
which is both Artinian and Noetherian. 
The simple objects are of the form 
$$
\iota_{!*}L[\dim_\C V]:={\rm Im}(\iota_{!}L\to\iota_{*}L)[\dim_\C V],
$$
where $V\hookrightarrow X$ is the immersion of locally closed subvariety of 
$X$ and $L$ is a local system on $V$.
\end{rem}

\begin{thm}$(${\bf Th\'eor\`eme 1.3.6 \cite{bbd:1}}$)$\\
The inclusion 
$ {}^pD^{\le 0}(\C_X) \hookrightarrow D_c^b(\C_X)$ 
$(resp. {}^pD^{\ge 0}(\C_X)\hookrightarrow D_c^b(\C_X))$ gives a 
right $($ resp.  left $)$ adjoint functor 
$ \tau_{\le 0}:D_c^b(\C_X) 
\longrightarrow  {}^pD^{\le 0}(\C_X) $,  
$ ($ resp. $ \tau_{\ge 0}:D_c^b(\C_X) \longrightarrow  {}^pD^{\le 0}(\C_X)$ $)$.

Moreover, 
$$
{}^p H^0:=\tau_{\ge 0}\tau_{\le 0}:D^b_c(\C_X)\to Perv(\C_X)
$$
is a cohomology functor, which  is called a {\em perverse cohomology functor}.

\end{thm}

By using the perverse cohomology functor, we 
can define the perverse direct images. 
\begin{defn}$(${\bf Perverse direct images functor}$)$\\
Let $f:X\to Y$ be a morphism of normal algebraic varieties.
$$
{}^pR^kf_*:Perv(\C_X)\to Perv(\C_Y),~~~K^\bullet\mapsto
{}^pR^kf_*K^\bullet:={}^pH^0(Rf_*K^\bullet[k]).
$$
\end{defn}

The following theorems are  the main results of 
\cite{bbd:1}. 
\begin{thm}\label{thm:decomp}$(${\bf Decomposition Theorem 
(Theor\'em\`e 6.2.5 \cite{bbd:1})}$)$\\
Let $f:X\to Y$ be a proper morphism of normal algebraic varieties and 
$K^\bullet\in Perv(\C_X)$ be a simple object of geometric origin.
Then 
$$
Rf_*K^\bullet\simeq\bigoplus_k 
{}^pR^kf_*K^\bullet[-k].
$$
\end{thm}

\begin{thm}\label{thm:relative}$(${\bf Relative hard Lefschetz theorem (Theor\'em\`e 5.4.10,  6.2.10 \cite{bbd:1})}$)$\\
Let $\omega$ be the first Chern class of the relative ample line bundle for the
projective morphism $f:X \to Y$.
Then for $k\ge 0$, we have
$$
\omega^k\wedge:{}^pR^{-k}f_*K^\bullet\simeq {}^pR^kf_*K^\bullet.
$$

\end{thm}

From Theorem \ref{thm:decomp} and \ref{thm:relative}, we can derive the following corollary. 
\begin{cor}\label{cor:fund}
Let $f:X \longrightarrow Y $ be a projective morphism between normal projective varieties. Moreover let $\omega_L$ and  $\omega_R$  
be the first Chern classes of a relatively 
ample invertible sheaf for $f:X \longrightarrow Y$ and 
an ample invertible sheaf of $Y$ respectively. 
Then the perverse Leray spectral sequence 
$$
{}^pE_2^{r, s}= {}^p H^r (Y, {}^p R^s f_* \C ) \Rightarrow {}^pH^{r+s}(X, \C) = IH^{r+s}(X, \C) 
$$
degenerates at $E_2$-term.  
Moreover  two 
relative hard Lefschetz actions for $f:X \longrightarrow Y $ and $ Y \longrightarrow \{ point \}$  
define actions on $E_2$-terms  
$$
\omega_L \wedge :{}^p E_2^{r,s} = {}^p H^r(Y, {}^p R^s f_* \C) \longrightarrow {}^p E_2^{r, s+2} = {}^p H^r(Y, {}^p R^{s+2} f_* \C)
$$
and 
$$
\omega_R \wedge:{}^p E_2^{r,s} ={}^p H^r(Y, {}^p R^s f_* \C) 
\longrightarrow {}^p E_2^{r+2,s} = {}^p H^{r+2}(Y, {}^p R^s f_* \C), 
$$
so that $(\omega_L)^s:{}^p E_2^{*, -s} \stackrel{\simeq}{\to} {}^p E_2^{*, s}$ and $(\omega_R)^{r}:{}^p E_2^{-r, *} 
\stackrel{\simeq}{\to} {}^p E_2^{r, *}$.  
These two actions commute to each other  and define an $(sl_2)_L \times (sl_2)_R $-action on the intersection cohomology ring $IH^{*}(X, \C)$. 
\end{cor}

\vspace{0.5cm}
\par\noindent
{\bf A digression for the representation of $sl_2$.}
\vspace{0.5cm}

We recall some fundamental facts on the representation of $sl_2$ and fix some notation.

It is well-known that the isomorphism class of the complex irreducible representations $V$ of $sl_2$ can be determined by their dimension $k$.
\begin{defn}
The irreducible representation of $sl_2$ of dimension $k$ is called the 
spin $\frac{k-1}{2}$ representation and it is denoted by 
$$
(\frac{k-1}{2}).  
$$
Note that for a non-negative 
half-integer $j \in \frac{1}{2} \Z_{\geq 0}$, 
the spin $j$ representation $(j)$ has dimension $2j +1$.  
\end{defn}  
Let $E, F, H$ be the usual generators of $sl_2$ with the relation 
$$
[E, F] = 2H, \ \ [H, E]= E, \ \ [H, F]= - F.
$$
For example, one may take a matrix representation 
$$
E =\left( \begin{array}{cc} 0 & 1 \\ 0 & 0 \end{array} \right),\quad F = \left( \begin{array}{cc} 0 & 0 \\ 1 & 0 \end{array} \right),  
\quad H = \frac{1}{2}\left( \begin{array}{cc} 1 & 0 \\ 0 & -1 \end{array} \right). 
$$
For the spin $j$ representation $(j)$, we can find an eigenvector $v \in V$ of $H$ with $F v = 0$.  Then one can show that 
$$
(j) = < v, E v, \cdots, E^{2j} v >_{\C}, \quad \mbox{and} \quad  E^{2j+1} v = 0, 
$$
$$
H (E^{k} v) = (- j + k) \cdot  v 
$$
In this case, the element $ E^{k} v $ has  the spin 
$ -j + k $, $ 0 \leq k \leq 2 j $.  

Let $E$ be an elliptic curve,  or a compact complex torus of dimension $1$ and let $\omega$ be a Lefschetz operator induced by a K\"ahler class.
  Then we have an isomorphism $H^2(E, \C) = \omega\cdot H^0(E, \C)$.     
The cohomology 
ring $H^*(E, \C)= H^0 \oplus H^1 \oplus H^2 \simeq \C \oplus \C^2 \oplus \C$ has the Lefschetz decomposition 
and it defines a
representation on $H^*(E, \C)$ of $sl_2$ as  
$$
(\frac{1}{2}) = H^0(E, \C)  \oplus \omega \cdot H^0(E, \C) \simeq \C \oplus \omega \C , \quad 2 \cdot (0) = H^1(E, \C) \simeq \C^2. 
$$  
Hence we have 
$$
H^*(E, \C) = (\frac{1}{2}) \oplus 2 (0). 
$$
In general, for a compact complex torus $A$ of complex dimension $g$, its complex cohomology 
ring has the Lefschetz decomposition induced by 
the Lefschetz operator $\omega$ which 
defines  a representation of $sl_2$ on $H^*(A, \C)$. It is obvious that 
as a representation of $sl_2$ 
$$
H^*(A,\C) = [ (\frac{1}{2}) \oplus 2 \cdot(0)]^{\otimes g}. 
$$
\begin{defn}
For each $g \geq 0$, we set 
$$
I_g = [ (\frac{1}{2}) \oplus 2 \cdot (0)]^{\otimes g}. 
$$
\end{defn}
\begin{rem}
For each irreducible representation $(j)$ of spin $j$, we can find integers $\alpha_r \in \Z$ so as 
\begin{equation}\label{eq:sl2}
(j) = \oplus_{r=0}^{2j} \alpha_r I_r.
\end{equation}
In fact, since $I_{2j}$ contains $(j) = H^0 \oplus \omega H^0 \oplus \cdots \oplus \omega^{2j} H^0 $ as an irreducible factor and the difference 
$$
I_{2j} - (j)
$$ 
is a sum of irreducible representations of spin $k$ where $ k < j$. By induction with respect to $j$, we see that 
 decomposition (\ref{eq:sl2}) holds. 
For example, we see that $(0) = I_0, (\frac{1}{2}) = I_1 - 2 I_0$,  $(1) = I_2 - 4I_1 + 3I_0, \cdots$.   
\end{rem}

Corollary \ref{cor:fund} says that 
for any projective morphism of projective normal varieties $ f : X \to Y $  we can define 
an action of $(sl_2)_L \times (sl_2)_R$ on the intersection cohomology ring $IH^*(X, \C)$. 
For any pairs $(j_1, j_2)$ of non-negative half-integers, let us set  $$
(j_1, j_2) = (j_1)_L \otimes (j_2)_R, 
$$
i.e., $(j_1, j_2) $ is an irreducible representation of $(sl_2)_L \times (sl_2)_L$ of bi-spin $(j_1, j_2)$.  Let us consider the irreducible decomposition of $IH^*(X, \C)$ 
defined by the relative Lefschetz action 
\begin{equation}
IH^*(X, \C) = \oplus \alpha_{j_1, j_2} (j_1, j_2).
\end{equation}
Moreover by the remark above, we can define the virtual decomposition of $IH^*(X, \C)$ as 
\begin{equation}
IH^*(X, \C) = \oplus_{h=0}^l I_h \otimes R_h
\end{equation}
where $R_j$ is a (virtual)  representation of $(sl_2)_R$.  

Summarizing the results, we obtain the following

\begin{thm}\label{thm:decomp2}
 Let $f:X \to Y$ be a projective morphism of normal complex projective varieties.  
Moreover let $\omega_L$ and $\omega_R$ be the first Chern classes of 
ample invertible sheaves on $X$ and $Y$, respectively.  Then the intersection cohomology ring $IH^*(X, \C)$ has a natural  $(sl_2)_L \times (sl_2)_R$ action induced by the relative Lefschetz action of $\omega_L$ and the Lefschetz action of $\omega_R$ on the $E_2$-terms of the perverse Leray spectral sequence for $f:X \to Y$.  The representation of $IH^*(X, \C)$ defines 
the irreducible decomposition and the virtual decomposition as 
\begin{equation}\label{eq:sl2decomp}
IH^*(X, \C) = \oplus_{j_1, j_2} \alpha_{j_1, j_2}(j_1, j_2) = \oplus_{h=0}^{l} I_h \otimes R_h. 
\end{equation}
\end{thm}

The following example shows that for a projective morphism $X \to Y$ the usual Leray spectral 
sequence does not detect the relative Lefschetz 
action and only the magic of the 
perverse sheaves can detect the natural 
$(sl_2)_L \times (sl_2)_R $-representation.  

\begin{exmp}$($Blowing up of $\P^2$$)$\\
Let $\pi:\widehat{\P}^2 \to \P^2$ be the blowing up at a point 
$ p \in \P^2$.

For $ \C_{\widehat{\P}^2}[2] \in  Perv (\C_{\widehat{\P}^2})$, the usual direct image $R^*\pi_{*}\C_{\widehat{\P}^2}[2]$ can be given by:

$$
R^j\pi_*(\C_{\widehat{\P}^2}[2]) = \left\{
\begin{array}{ll}
\C_{\P^2}  \quad & \mbox{for} \quad  j = -2 \\
\C_{p} \quad  & \mbox{for} \quad  j = 0 \\
0  \quad  & \mbox{otherwise}
\end{array} \right.
$$
The $E_2$-term of the ordinary Leray spectral 
sequence is given by 
$$
E^{i, j} =  
H^{i}(\P^2, R^j \pi_*( \C_{\widehat{\P}^2}[2])), 
$$ 
and more explicitly each term  is  given in 
the  following table.

\begin{table}[h]
\begin{center}
\caption{}
\label{tab:ordinary}
\begin{tabular}{|cccc} 
 $E^{4, -2}$ & $E^{4, -1}$ & $ E^{4, 0}$ &   \\
$E^{3, -2}$ & $E^{3, -1}$ &  $E^{3, 0}$ &  \\
$E^{2, -2}$ & $E^{2, -1}$ &  $E^{2, 0}$ &  \\
 $E^{1, -2}$ & $E^{1, -1}$ & $ E^{1, 0}$ &  \\
 $E^{0, -2}$ & $E^{0, -1}$ &  $E^{0, 0}$ &   \\ \hline 
\end{tabular} 
\quad 
$=$ \quad \begin{tabular}{|ccc}
    $\C$ & $0$ & $0$  \\
    $0$ & $0$ & $0$ \\
    $\C$ & $0$ & $0$ \\
    $0$  & $0$ & $0$  \\
    $\C $ & $0$ & $\C$ \\  \hline 
 \end{tabular} 
\end{center}
\end{table}
\begin{equation}
\begin{array}{lcl}
\H^{-2}(\widehat{\P^2}, \C[2])\simeq H^0(\widehat{\P^2}, \C)  & \simeq  & E^{0, -2} = \C \\
\H^{0}(\widehat{\P^2}, \C[2])\simeq H^2(\widehat{\P^2}, \C)  & \simeq & E^{2, -2} \oplus E^{0,0} \simeq \C \oplus \C \\
\H^{2}(\widehat{\P^2}, \C[2]) \simeq H^4(\widehat{\P^2}, \C)  & \simeq & E^{4, -2}  \simeq  \C \\
\end{array}
\end{equation}

From this diagram, one can not detect a  natural 
action of $(sl_2)_L \times (sl_2)_R$ on the cohomology 
ring $H^*(\widehat{\P^2}, \C)$ because the Table 1 is not 
symmetric.  
\
On the other hand, the perverse direct image can be given 
as 
$$
{}^pR^j\pi_* \C_{\widehat{\P}^2}[2] =
\left\{ \begin{array}{ll} \C_{\P^2}[2] \oplus \C_{p} & \quad \mbox{for} \quad j = 0  \\
0  & \quad \mbox{otherwise} \end{array} \right.  
$$
that is, the complex is the direct sum of $\C_{\P^2}[2] \oplus \C_{p}$ and 
 concentrated on  the degree $0$  and
moreover the relative Lefschetz action 
(= the left action of   $ \omega_L $ )
on $ {}^p R^* \pi_* \C_{\widehat{\P}^2}[2]= \C_{\P^2}[2] \oplus \C_{p}$ is trivial.     
Hence we have the decomposition 
$$
\H^*(\widehat{\P^2}, \C_{\widehat{\P}^2}[2]) \simeq \H^*(\P^2, {}^p R^*\pi_{*}\C_{\widehat{\P}^2}[2]) = \H^* (\P^2, \C_{\P^2}[2]) \oplus 
\H^*(p, \C_p)
$$
where each decomposition factor has the right $sl_2$-action 
as
$$
\begin{array}{cccccc}
\H^{-2}(\P^2, \C_{\P^2}[2]) &  \stackrel{\omega_R \wedge}{\longrightarrow} & \H^{0}(\P^2, \C_{\P^2}[2]) &  \stackrel{\omega_R \wedge}{\longrightarrow} & \H^2
(\P^2, \C_{\P^2}[2]) & \mbox{(spin $1$) } \\
 || &  & ||  & & ||  \\
 \C &  \stackrel{\omega_R \wedge}{\longrightarrow} & \C \cdot \omega_R  & \stackrel{\omega_R \wedge}{\longrightarrow} &  \C \cdot \omega_R^2   
\end{array}
$$ 
and 
$$
\begin{array}{cccccc}
\H^{-2}(p, \C_p)  & \stackrel{\omega_R \wedge}{\longrightarrow} & \H^{0}(p, \C_p) & \stackrel{\omega_R\wedge}{\longrightarrow} & \H^2
(p, \C_p) & \mbox{(spin $0$)}  \\
 || &  & ||  & & || &  \\
 0  &  \stackrel{\omega_R \wedge}{\longrightarrow} & \C   &\ \stackrel{\omega_R \wedge}{\longrightarrow} &  0 . & 
\end{array}
$$
Noting that the action of $\omega_L$ is trivial, we see that
the following isomorphisms of  $(sl_2)_L \times (sl_2)_R$ representations 
$$
\H^{*}(\P^2, \C[2]) \simeq (0, 1) \simeq (0)_L \otimes (1)_R, \quad 
\H^{*}(p, \C_p) \simeq (0, 0) \simeq  (0)_L \otimes (0)_R. 
$$
As a result, the $(sl_2)_L \times (sl_2)_R$ 
decomposition  of  
$H^*(\widehat{\P^2},\C)$ is given by 
\begin{eqnarray*}
H^*(\widehat{\P^2},\C)& = & \H^{*}(\P^2, \C[2]) 
\oplus \H^{*}(p, \C_p)  \\
&   \simeq &   (0)_L\otimes [(1)_R\oplus(0)_R]= I_0 \otimes [(1)_R\oplus(0)_R]  
\end{eqnarray*}

For reader's convenience, we put a table of  
the $E_2$-terms of the perverse Leray spectral sequence. 
One  can  compare these with those of the ordinary one in Table \ref{tab:ordinary}:
\begin{table}[h]
\caption{}    
\begin{center}
\begin{tabular}{c|cccc} 
&${}^p E^{2, -1}$ & $ {}^pE^{2, 0}$ & $ {}^p E^{2, 1}$ &   \\
&${}^pE^{1, -1}$ & ${}^pE^{1, 0}$ &  ${}^pE^{1, 1}$ &  \\
 $\uparrow  \omega_R \wedge $ &${}^pE^{0, -1}$ &  ${}^pE^{0, 0}$ &  ${}^pE^{0, 1}$ &  \\
&${}^pE^{-1, -1}$ & ${}^pE^{1, 0}$ & $ {}^pE^{-1, 1}$ &  \\
&${}^pE^{-2, -1}$ & ${}^pE^{-2, 0}$ &  ${}^pE^{-2, 1}$ &   \\ \hline 
  &  & $ \omega_L \wedge  \to$ &  &   
\end{tabular} 
\quad 
$=$ \quad \begin{tabular}{c|cccc}
   & $0$ & $\C \oplus 0 $ & $0$ &   \\
   &  $0$ & $0$ & $0$ &  \\
$ \uparrow  \omega_R \wedge $  &  $0$ & $\C \oplus \C$ & $0$ & \\
   & $0$  & $0$ & $0$  & \\
   &  $0 $ & $\C \oplus 0 $ & $0$  & \\  \hline 
   &       &     $  \omega_L \wedge \to$          &    & 
 \end{tabular} 
\end{center}
\end{table}

\end{exmp}
%%%%%%%%%%%%%%%%%%%%%%%%%%%%%%%%%%%%%%%%%%%%%%%%%%%%%%%%%%%%%%%%%%%%%%%%%%%%%%
%%%%%%%%%%%%%%%%%%%%%%%%%%%%%%%%%%%%%%%%%%%%%%%%%%%%%%%%%%%%%%%%%%%%%%%%%%%%%%
\section{D-brane moduli spaces}
~~~

\vspace{0.2cm}
What is the mathematical definition of  ``D-branes wrapped around a cycle of dimension one" and the moduli space $ M $ of ``D-branes"?
Usually one may think them as cycles with flat $U(1)$-bundles (or equivalently holomorphic line bundles of degree $0$) and the 
moduli of them. 

This translation is sufficient in many cases, but since the cycles may have singularities, it is more useful for our purpose to regard D-branes as stable sheaves (Narasimhan--Seshadri theorem, 
Kobayashi--Hitchin correspondence). 

Moreover, it is rather subtle to define the moduli of supports of 
sheaves and the natural support map 
\begin{equation}\label{eq:sup}
\pi:M \lra S
\end{equation}
from the moduli space $M$ of sheaves to the moduli space  $S$ of cycles.  (cf. \cite{lp:1}).  

For example, let us consider the following situations:
Let $X$ be a Calabi-Yau threefold and let $C \subset X$ be a 
smooth irreducible curve.  We can consider the following two cases:

\begin{enumerate}
\item $n$ copies of D-branes wrapped around the cycle $C$ once.
\item Large single D-brane wrapped around the cycle $C$ $n$-times.
\end{enumerate}

Mathematically, the first one corresponds to a sheaf 
of rank $n$ on $C$ 
and the second one corresponds to a sheaf 
of rank $1$ on non-reduced scheme with 
the same topological space $C$ (but multiplicity along $C$ is $n$). 
Sometimes the above two objects have the same Hilbert polynomial 
and hence they define points in  the same moduli space and the first 
one may be deformed to the second one algebraically.  

Therefore, in order to make the support map (\ref{eq:sup}) 
a morphism of algebraic schemes, 
one has to define  natural multiplicities  of 
irreducible components of the support of the corresponding sheaf.  
    Since there seems to be no natural way to put the 
scheme structure on the supports of pure sheaves when the support 
has codimension greater than one, we consider the supports of sheaves
with multiplicities as the algebraic cycles in the total space $X$.  
Hence the moduli space $S$ can be considered as a subset of Chow varieties  
of $X$ and the support map (\ref{eq:sup}) can be considered as  
 a generalization of 
the morphism from Hilbert scheme to Chow varieties.

Let us first recall the necessary background in the theory of  moduli spaces of sheaves (cf. \cite{hl:1}).  
  Let $Z$ be a Noetherian scheme and $\E$ be a 
coherent sheaf on $Z$.  
\begin{defn} The support of $\E$ is the close subset ${\rm Supp}( \E) = \{ z \in Z |  \E_{z} \not= 0 \} $. ${\rm Supp}(\E)$ becomes a closed reduced subscheme of $X$.  
 Its dimension is called the dimension of the sheaf $\E$ and is denoted by $\dim (\E)$.  
\end{defn}

\begin{defn} 
A coherent sheaf $\E$ on a scheme $Z$ is pure of dimension $k$ 
if $\dim_\C{\rm Supp}({\mathcal F})=k$ for any nontrivial coherent subsheaf 
${\mathcal F}\subset\E$.
\end{defn}
\begin{defn} Let $Z$ be a projective scheme over $\C$ and $L$ be an ample line bundle on $Z$, let  $\E$ be a coherent sheaf which is pure of dimension $d$ on $Z$.   Let 
$$
P(\E,m):=\chi(Z,\E(m))=\sum_{i=0}^d\alpha_i(\E)\frac{m^i}{i!}
$$
be the Hilbert polynomial of $\E$.  (Here $\E(m) := \E \otimes L^m$.) 
Then $p(\E,m):=P(\E,m)/\alpha_d(\E)$ is called a reduced Hilbert polynomial 
of $\E$.
\end{defn}

\begin{defn}$(${\bf Stability}$)$  Let $\E$ be a coherent sheaf which is pure of dimension $d$ on a projective 
scheme $Z$.
$\E$ is stable (resp. semistable) 
if for any proper subsheaf ${\mathcal F}$, 
$$
p({\mathcal F},m)< p(\E,m),~~~{\rm for}~~m>>0.
$$
$$
(\mbox{resp.} \quad p({\mathcal F},m) \le p(\E,m),~~~{\rm for}~~m>>0).
$$
\end{defn}

Let $X$  be a smooth projective scheme over $\C$ and 
let $\E$ be a coherent sheaf on $X$  which is of pure of dimension 1. 
Let  ${\rm Supp} (\E)$ be the support of $\E$ and let  
 $Y_1, \cdots,  Y_l $ be the irreducible components 
 of ${\rm Supp}(\E)$, let $v_i$ be the generic point of $Y_i$.  
 Then the stalk $\E_{v_i} = \E \otimes_{\cO_{X}} \cO_{X, v_i}$ 
 is an Artinian module of finite length $ length (\E_{v_i})$. 
 We define an algebraic cycle  $s(\E)$ by 
 \begin{equation}\label{eq:cycle}
 s(\E) :=  \sum_{i=1}^l length (\E_{v_i})  \cdot Y_i. 
 \end{equation}
Moreover the homology class of $s(\E)$ will be denoted by 
$[s(\E)] \in H_2(X, \Z)$.

We can define the following moduli spaces of semistable sheaves by 
the Simpson's construction (see, for example \cite{hl:1}):

\begin{defn} Let $X$ be a smooth projective 
 Calabi-Yau 3-fold defined over $\C$ and $L$  an ample line bundle on $X$.  For a positive integer $d$, let   $M_{d}(X) $ be  the moduli space of semistable sheaves 
$\E$ of pure dimension 1  on $X$ with Hilbert polynomial 
\begin{equation}
P(\E,m)= d \cdot m +1
\end{equation}
It is known that $M_d(X)  $ is a projective scheme over $\C$ (Theorem 4.3.4 \cite{hl:1}).
\footnote{From the view point of physics, we may consider 
 $M_\beta= M_{\beta}(X)$ 
as the moduli space of D-branes by the physical 
discussion that the degeneracy of BPS states should be 
independent of the $U(1)$ flux (the degree of sheaves).} 
\end{defn}

Let $\E \in M_{d}(X)$.  It is easy to see that the degree $d$ of $\E$ is given by the intersection number
$$
d = L \cdot s(\E) = [L] \cdot [s(\E)]
$$
of $L$ with the support cycle $s(\E)$. (Here 
$[L]$ denote the homology class of the divisor associated to $L$).    
Let $Chow_d(X)$ denote the Chow variety parameterizing  
algebraic 1-cycles $Y$ on $X$ with $ L \cdot Y = d$.  It is known that $Chow_d(X)$ is a projective scheme over $\C$.  

In the same way as (\S 5, Ch.\ 5, \cite{M}), we can show that  the natural map 
\begin{equation}\label{eq:sheadtochow}
\begin{array}{cccc}
\pi_{d}:& M_{d}(X) & \lra &  Chow_{d}(X). \\
        &  \E & \mapsto & s(\E)  
\end{array}
\end{equation}
becomes a morphism of projective schemes.  
 
For a  homology cycle $\beta \in H_2(X, \Z)$ with $ L \cdot \beta = d$, we can define the closed subscheme
\begin{equation}
M_{\beta}(X) := \left\{ \begin{array}{c|l}  
 & \mbox{a semistable sheaf on $X$, pure of  } \\
\E  & \mbox{dimension $1$ with} \  P(\E, m) = d m + 1 \\
  & and   \  [s(\E)] = \beta  \end{array}     \right\}/ \mbox{isom.} \subset M_d(X). 
\end{equation} 
 \begin{equation}
 Chow_{\beta}(X) := \{  
\gamma \in Chow_d (X),  \ [\gamma] = \beta \} \subset Chow_d(X). 
\end{equation} 
Then we can define  a natural morphism
\begin{equation}
\pi_{\beta}:M_{\beta}(X) \lra Chow_{\beta}(X)
\end{equation}
by $\pi_{\beta}(\E) = s(\E)$.  
Taking the normalization of 
$\tilde{M}_{\beta}(X) \lra M_{\beta}(X)$ and 
setting 
$$
S_{\beta}(X) = \mbox{the normalization of the image 
$ \pi_{\beta}(M_{\beta}(X))$ }, 
$$ 
from the univesal property of the normalization, 
we obtain  a natural 
surjective morphism between normal projective varieties over $\C$:
\begin{equation}\label{eq:map}\boxed{ \ \ 
\pi_{\beta}:\tilde{M}_{\beta}(X) \lra S_{\beta}(X). \ \ }
\end{equation}

Our proposal for  definition of BPS invariants 
which may be  consistent  with  
the Gopakummar--Vafa conjecture \ref{conj:gvc} is that 
the moduli $\tilde{M}_{\beta}(X) $ of sheaves on $X$ with the 
homology class of support cycle $\beta$  should be the natural 
moduli of $D$-branes wrapping  around a support cycle $\beta$.  
We state this as our conjecture.

\begin{conj}
The morphism \eqref{eq:map})$\pi_{\beta}:\tilde{M}_{\beta}(X)  
\lra S_{\beta}(X) $ is the natural 
morphism from $D$-brane moduli space $\tilde{M}_{\beta}(X) $ to the 
moduli space $S_{\beta}(X)$ of support curves whose homology class is 
$\beta$. 
\end{conj}

Now as suggested in \cite{gv:1}, 
applying  Theorem  \ref{thm:decomp2} to   the morphism 
\eqref{eq:map}, we 
obtain the following theorem and definition of BPS invariants.  
\begin{thm} \label{thm:bps}
Let $\pi_{\beta}:\tilde{M}_{\beta}(X) \longrightarrow S_{\beta}(X) $ be the projective morphism defined as in \eqref{eq:map} and fix ample line bundles $L_1$ on $\tilde{M}_{\beta}(X)$ and $L_2$ on 
$S_{\beta}(X)$ respectively.  
  Then there exists an $(sl_2)_L \times (sl_2)_R$-action on 
$IH^*(\tilde{M}_{\beta}(X))$ defined by the relative Lefschetz operator $\omega_L$ and by the Lefschetz operator $\omega_R$ of the base.
The $(sl_2)_L \times (sl_2)_R$-action gives the decomposition of 
$IH^*(\tilde{M}_{\beta}(X))$
\begin{equation}\label{eq:decomp-sl2}
IH^*(\tilde{M}_\beta(X))= \bigoplus_{h\ge 0} I_h \otimes R_h(\beta)
= 
\bigoplus_{h\ge 0}
\left[(\frac{1}{2})_L\oplus 2(0)_L\right]^{\otimes h}\otimes R_{h}(\beta).
\end{equation}
where we denote by $(j)_L$ the spin-$j$ 
representation of the relative Lefschetz 
$(sl_2)_L$-action and by 
$R_h(\beta)$ a (virtual) representation of the 
$(sl_2)_R$-action.
\end{thm}

\begin{defn}\label{def:bps}
By using the decomposition (\ref{eq:decomp-sl2}), 
we can define integers  $n_h(\beta)$, which will be 
called {\em BPS invariant},   by the following formula:
\begin{equation}\label{eq:bps}
 \boxed{ \  n_{h}(\beta):=Tr_{R_{h}(\beta)}(-1)^{2H_R}. \ }
\end{equation}
\end{defn}

\begin{conj}\label{conj:gv2}
Integers $n_h(\beta)$ defined in \eqref{eq:bps} 
 should be deformation invariants satisfying the 
conjecture  $\ref{conj:gvc}$ of Gopakumar--Vafa. 
In particular, $n_0(\beta)$ should be the holomorphic Casson invariants defined by Thomas \cite{th:1}.
\end{conj}

Since neither  $\tilde{M}_\beta(X)$  or  the morphism $\pi_\beta$ may 
not be  smooth in general, 
we cannot prove the existence of such an action on 
$H^*(\tilde{M}_\beta(X), \C )$ by the usual Leray's
 spectral sequence. 
However, the ``perverse" Leray spectral sequence tells us the origin of the $(sl_2)_L\times(sl_2)_R$-action on intersection cohomology 
$IH^*(\tilde{M}_\beta(X))$. 
(Note that  if $\tilde{M}_\beta(X)$ is smooth,   $IH^*(\tilde{M}_\beta(X))=H^*(\tilde{M}_\beta(X))$).

%%%%%%%%%%%%%%%%%%%%%%%%%%%%%%%%%%%%%
%%%%%%%%
%
%%%%%%%%%%%%%%%%%%%%%%%%%%%%%%%%%%%%%%%%%%%%%%%%%%%%%%%%%%%%%%%%%%%%%%%%%%%%%%
%%%%%%%%%%%%%%%%%%%%%%%%%%%%%%%%%%%%%%%%%%%%%%%%%%%%%%%%%%%%%%%%%%%%%%%%%%%%%%

\vspace{0.5cm}
\section{Evidences}

In Section 3, we gave a mathematical definition of BPS invariants $n_h(\beta)$ (cf. Definition \ref{def:bps}).  
In this section, we will present  several pieces of  evidences 
supporting  Conjecture \ref{conj:gv2} or  Conjecture \ref{conj:gvc} by using Definition \ref{def:bps}.  
Rigid rational curves,  super-rigid elliptic curves and also some curves in rational 
elliptic surfaces will be considered.

\subsection{Super-rigid curves in a Calabi--Yau $3$-fold and conjectural local BPS invariants}

\quad 

Let $C \subset X$ be a smooth irreducible curve in 
a Calabi-Yau 3-fold with the homology class 
$\beta = [C] \in H_2(X, \Z)$.  The curve $C \subset X$ is {\em rigid} 
if $H^0(C, N) = 0$ where $N = N_{C/X}$ be the normal bundle of $C$ in $X$.  Moreover $C \subset X$ is called {\em super-rigid} if, 
for all non-constant maps of nonsingular curves $\mu:C' \longrightarrow C$, 
\begin{equation}\label{eq:s-r}
H^0(C', \mu^*(N)) = 0. 
\end{equation}
For a super-rigid curve $C$ of genus $ 0$ or $1$, one can 
define the local contributions  of $C$ to the Gromov-Witten invariants $N_g (n [C]) = N_g (n \beta)$ (\cite{pa:1}, \cite{bp:1}).  Denote corresponding local Gromov--Witten invariant by $N_g (n \cdot C)$. 
Define the {\em (conjectural) local BPS invariants} $n_g^{conj}( d \cdot C) \in \Q $ by  (conjectural) Gopakummar--Vafa 
formula \eqref{eq:gv} 
\begin{equation}
\sum_{g\ge 0, n \ge 0 }N_g(n \cdot C )q^n \lambda^{2g-2}=
\sum_{k>0, g \ge 0, n \ge 0 } n_g^{conj}(n \cdot C)\frac{1}{k}
\left(2\sin(\frac{k\lambda}{2})\right)^{2g-2}q^{kn}.
\end{equation}
(Matching the coefficients of the two series yeilds equation 
determining $n_g^{conj}(n \cdot C)$ recursively in terms of $N_g(n \cdot C )$ (cf. Proposition 2.1 \cite{bp:1}). Note that in the notation of \cite{bp:1} one has $N_{g+h} (n \cdot C) = N_n^{h} (g)$ and $n_{g+h}^{conj}( n \cdot C) = n_n^h(g)$ where 
$g=$ genus of $C$).

For the case of $C = \P^1$, Faber and Pandharipande proved the following theorem for the generating function of local Gromov--Witten 
invariants $N_g(n \cdot \P^1)$:
\begin{thm}\label{thm:srr}
$($\cite{fp:1}$)$
$$
\sum_{g\ge 0,n\ge 1}N_g(n \cdot \P^1)q^n\lambda^{2g-2}=
\sum_{k\ge 1}\frac{1}{k}\left(2\sin(\frac{k\lambda}{2})\right)^{-2}q^k.
$$
\end{thm}

From this formula for Gromov--Witten invariants, 
the conjectural local BPS invariants $n_g^{conj}(n \cdot \P^1)$ 
can be given by (cf. \cite{bp:1}) 
\begin{equation}\label{eq:bps-rr}
n_{g}^{conj}(n \cdot \P^1) = 
\left\{\begin{array}{l}
1 \quad \mbox{for  $g=0$ and $n = 1$} \\
0 \quad \mbox{otherwise. }
\end{array} \right. 
\end{equation}
 
Next let  $E \subset X$ be a super-rigid elliptic curve. 
 Pandharipande \cite{pa:1}  showed the following:
 
\begin{thm}$ \label{thm:sre}
($\cite{pa:1}$)$ \\
Let $E \subset X$ be a super-rigid elliptic curve.   One has 
$$
N_1( n \cdot E) = \frac{\sigma(n)}{n} = \sum_{i|n} \frac{1}{i}.  
$$
Moreover for all $g > 1,   n >0 $, one has 
$$
N_g(n \cdot E) = 0. 
$$
Therefore  
the LHS of conjectural formula \eqref{eq:gv} are given by 
$$
\sum_{n\ge 1}N_1(n \cdot E)q^n=-\sum_{n\ge 1}\log(1-q^n)
=\sum_{n \ge 1,k\ge 1}\frac{1}{k}q^{kn}.
$$
\end{thm}

Again the conjecture \eqref{eq:gv} with Theorem \ref{thm:sre} reads 

\begin{equation}\label{eq:se}
n_g^{conj}( n \cdot E) = 
\left\{\begin{array}{l}
1 \quad \mbox{for  $g=1$ and all $n  \geq 1$} \\
0 \quad \mbox{otherwise. }
\end{array} \right. 
\end{equation}

\subsection{Calculations of BPS invariants $n_g( d \cdot C)$ via 
moduli of sheaves} 

\quad
\vspace{0.2cm}

   Let $X$ be a Calabi-Yau 3-fold and fix an ample line bundle $\O_X(1)$ on $X$ with $\omega := c_1(\O_X(1))$.  
Let $C \subset X$ be a super-rigid rational or elliptic curve of degree 
$$
d = \O_X(1) \cdot C = \int_{[C]} \omega.
$$  
Let $[C] \in H_2(X, \Z)$ denote the 
homology class of $C \subset X$.  In Definition \ref{def:bps}, for a 
non-negative integer $g$, 
we define the BPS invariant $n_g( n \cdot [C])$ by \eqref{eq:bps}
with respect to the 
surjective projective morphism defined in \eqref{eq:map} 
\begin{equation}
\pi_{n\cdot [C]}: \tilde{M}_{n \cdot [C]} \lra  S_{n \cdot [C]}. 
\end{equation}

Our next aim is to calculate  (local) 
BPS invariants $n_g (d \cdot C) $ defined by  (\ref{eq:bps}) and 
to compare $n_g (d \cdot C) $ with $n_g^{conj}(d \cdot C)$.  

In order to consider the local BPS invariant $n_g (d \cdot C) $, 
let us consider the  subset (or more explicitly the subfunctor) 
 of $M_{n [C]}(X)$ defined by 
\begin{equation}
M_{n \cdot C}(X)  :=  \left\{ \quad   \E \quad 
\begin{array}{|l} \mbox{a stable coherent $\O_X$-sheaf, pure of dimension 1} \\
\mbox{with $supp(\E) = C$ and $P(\E, m) =  n d m   +1$. }
\end{array}
 \right\} /\mbox{isom.}
\end{equation}
Note that $\E \in M_{n \cdot C}  (X)  $ implies that the cycle 
theoretic support $s(\E)$ of $\E$ is given by $n \cdot C$.

In Proposition \ref{prop:srrational} and \ref{prop:srelliptic}, we will show that  $M_{n \cdot C}(X) $ is a smooth irreducible component 
of $M_{ n \cdot [C]} (X) $ or  the empty set.  
Therefore if $M_{ n \cdot [C]} (X) $ is not empty set, 
the image of the map 
$$
\pi:M_{n \cdot C}(X)  \lra  S_{n \cdot [C]}(X).  \quad 
\pi(\E) = s(\E) = n \cdot C
$$
consists of  just one point $\{ s(\E) = n \cdot C \} \simeq 
\Spec \C$.  We remark that since $M_{n \cdot C}(X) $ is smooth, 
we will not distingusih   $M_{n \cdot C}(X)$ and its normalization  
$\tilde{M}_{n \cdot C}(X)$.

Let us set $M_{n \cdot C} = M_{n \cdot C} (X)$ and consider 
the natural morphism $\pi:M_{n \cdot C} \lra  \Spec \C$.   
By using the decomposition of the intersection cohomology 
ring $IH^*(M_{n \cdot C}, \C)$ with respect to 
the morphism $\pi$ (cf.\ Theorem \ref{thm:bps}), we obtain the decomposition into  
$(sl_2)_L \times (sl_2)_R$-representations as 
$$
IH^*(M_{n \cdot C}, \C)= \oplus_{h \geq 0} I_h \otimes R_h(n \cdot C).
$$
Now the local BPS invariants $n_{h}(d \cdot C)$ are given by 
(cf.\ Theorem \ref{thm:bps}) 
$$
n_{h}(n  \cdot C) = Tr_{R_h(n \cdot C)} (-1)^{2H_R}. 
$$

Now let us consider the moduli space $M_{n \cdot C}(X)$ for a 
super-rigid rational or elliptic curve $C \subset X$.  

\begin{prop}\label{prop:srrational}
Let $X$ be a Calabi-Yau 3-fold, $C \subset X$  a smooth  rigid rational curve (i.e., $C = \P^1$ and $N = N_{C/X} = \cO_C(-1) \oplus \cO_C(-1)$).  
Let $\omega$ be the class of an ample line bundle $L$ 
 on $X$ and set 
$d = L \cdot C = \int_{\beta} \omega$.  Let $\E$ be an element of 
the moduli space $M_{n \cdot C}(X) = M_{n \cdot \P^1}(X)$, that is, 
 an isomorphism class of a  
stable $\O_X$-coherent sheaf with $supp(\E) = C$ and $P(\E, m) =  
n d m + 1$ for some positive integer $n$.  Then we have 
$$
n = 1, \quad \mbox{and} \quad  \E = \O_C.  
$$
Moreover we have isomorphisms of schemes 
\begin{equation}
M_{n \cdot \P^1} \simeq 
\left\{ \begin{array}{ll}
 \{ \O_C  \} \simeq \Spec \C  \quad  \mbox{(one point)} & \quad \mbox{if $n =1$} \\
\quad  \emptyset  & \quad \mbox {otherwise. } 
 \end{array} \right.
\end{equation}
\end{prop}

\vspace{0,5cm}
\noindent
{\it Proof. \footnote{We thank Kota Yoshioka for suggesting us the following proof and also the proof of Proposition \ref{prop:srelliptic}. } } 
Let $\E$ be a stable coherent $\O_X$ sheaf with $supp(\E) = C$ and 
 $P(\E, m) = n d m + 1$.  
 Since $\dim H^0(X, \E) \geq P(\E, 0) = 1$, there exists a non-trivial section $s$ of $\E$ which defines a morphism of $\O_X$-sheaves 
$$
 s: \O_X \longrightarrow \E.
$$
Let $J$ be the kernel of the morphism $s$ and $I$ the ideal sheaf of $C \subset X$.  
Since $supp (\E) = C$ and $C$ is an irreducible reduced curve, we see that $J \subset I$.  On the other hand, the morphism $s$ defines an 
injection
\begin{equation}\label{eq:inj}
\varphi_s: \O_X/J \hookrightarrow \E.
\end{equation}
First assume that $J = I$.  Then we have an injection $\varphi_s: \O_C \hookrightarrow \E$.  If $ \O_C \subsetneqq \E $, the stability condition implies that 
$$
p(\O_C, m )  = \frac{1}{d} \chi(\O_C \otimes L^m) < p(\E, m) = \frac{1}{nd} \chi(\E \otimes L^m )
$$
or 
$$
m + \frac{1}{d} < m + \frac{1}{nd}.
$$
Since $d > 0$ and $n  \geq 1$, this gives the contradiction.  
Hence this implies that if $ I = J $ $ \E = \O_C$ and $n =1$. 

Next let us consider the case $ J \subsetneqq I$. We have an integer $k \geq 1$ such that $I^{k+1} \subseteqq J$ and $I^k \not\subseteqq J $.
Then  $(J + I^k)/J$ is a non-trivial subsheaf of $\O_X/J$ and 
hence $(J + I^k)/J \subseteq \E$. 
By using the isomorphism $(J+I^k)/J \simeq I^k/J \cap I^{k}$ and $
I^{k+1} \subseteq J \cap I^{k} $, 
we see that there exists a surjection 
$$
I^k/I^{k+1} \simeq S^k( \O_C(1)^{\oplus 2}) \twoheadrightarrow I^k/J \cap I^{k}. 
$$
(Note that $I/I^2 \simeq N^{\vee} \simeq \O_C(1)^{\oplus 2}$.  
Here $\O_C(1)$ is the ample generator of $\Pic(C) = \Z$. ) 
Consequently, one obtains the non-trivial morphism 
$$
\O_C(k ) \longrightarrow \E, 
$$
which again contradicts the stability of $\E$.  Thus $\E \simeq O_C$ and $n=1$. 

We have proved that the set of $\C$-valued point $M_{n \cdot \BP^1}$ 
consists of $\O_C$ if $n =1$ and is empty if $ n > 1$.  
In order to see that $M_{1 \cdot \BP^1}$ is smooth, it suffices to show that  the Zariski tangent space $T_{[\O_C]}$ of $M_{1\cdot \BP^1}$ at  $\O_E$ is $0$.  By general theory, we have the isomorphism 
\begin{equation}\label{eq:zariski}
T_{[\O_C]} \simeq \Ext^1_{\cO_X} (\O_C, \O_C). 
\end{equation}
Hence it suffices to show that 
\begin{equation}\label{eq:ext-van-r}
\Ext^1_{\cO_X} (\O_C, \O_C) \simeq \Ext^1_{\O_C}(\O_C, \O_C) \simeq H^1(C, \O_C)  \simeq H^1(\BP^1, \O_{\BP^1}) \simeq 0.  
\end{equation}

Recall the exact sequence induced by local to global spectral sequence

\begin{equation}\label{eq:l-g-r}
\begin{array}{l}
0 \rightarrow 
 H^1( \cHom_{\O_X}(\O_C, \O_C))  \rightarrow  
\mbox{\rm Ext}^1_{\O_X}(\O_C, \O_C)   \rightarrow  H^0( \cExt^1_{\O_X}(\O_C, \O_C))  \\ 
   \quad \quad  \rightarrow 
  H^2( \cHom_{\O_X}(\O_C, \O_C)) 
  \end{array}
\end{equation}
From the exact sequence 
$$
 0 \lra I \lra \O_X \lra \O_C \lra 0,  
$$
and the isomorphism  $ \cHom_{\O_X}(\O_C, \O_C) \stackrel{\simeq}{\lra}  \cHom_{\cO_X}(\O_X, \O_C) \simeq \cO_C$, we obtain the isomorphism  
$$
\cHom_{\O_X}(I, \O_C) \simeq  \cExt^1_{\O_X}(\O_C, \O_C). 
$$
Since $\cHom_{\O_X}(I, \O_C) \simeq \cHom_{\O_C}(I/I^2, \O_E)$, we finally obtain the isomorphism 
\begin{equation}
\cExt^1_{\O_X}(\O_C, \O_C) \simeq \cHom_{\O_C}(I/I^2, \O_C) \simeq N 
\end{equation}
where $N = N_{C/X}$ is the normal sheaf for $C \subset X$. 
Therefore we have 
$$
H^0( \cExt^1_{\O_X}(\O_C, \O_C)) \simeq H^0(\BP^1, N) \simeq H^0( \BP^1, \O_{\BP^1}(-1)^{\oplus 2}) \simeq 0. 
$$
Hence the sequence \eqref{eq:l-g-r} implies  \eqref{eq:ext-van-r}.

\qed

\vspace{0.5cm}
\noindent
Let $ E \subset X$ be a super-rigid elliptic curve. Then it is 
known that the normal bundle $N$ of $E$ in $X$ is isomorphic to 
$$
N \simeq L \oplus L^{-1}
$$  
where $L$ is a non-torsion element of 
the Picard group of $E$, the group of line bundle on $E$ of 
degree $0$ .  

Let $ M_{n, 1}(E)$ denote the moduli space of stable $\O_E$-locally free sheaves  of rank $n$ of degree $1$.  Fix an ample line bundle $L$ on $X$ and denote it by $\omega = c_1(L)$ and   set $d = L \cdot E = \int_{[E]} \omega$.  We denote by  $M_{n \cdot E}$ the relevant moduli space  of stable $\O_X$ sheaves $\E$  with $supp(\E) = E $ and $P(\E, m) = \chi(\E \otimes L^m) = n d m + 1$.  Let $\iota: E \hookrightarrow X$ be the natural inclusion.  Then we have a natural push-forward morphism
\begin{equation}\label{eq:srelliptic}
\begin{array}{ccc}
 \iota_{*}: M_{n, 1}(E) & \longrightarrow &  M_{n \cdot E} \\
        &   &   \\
  \quad \F & \mapsto &  \iota_*(\F) . 
\end{array}
\end{equation}

Atiyah ([Theorem 7, \cite{A}])  showed  that $ M_{n, 1}(E)$ is isomorphic to $M_{1, 1}(E)$ via the map $\E \mapsto \det \E$. Hence 
fixing a line bundle $A$ of degree $1$ and identifying  $M_{1, 1}(E)$ with  $ M_{1,0}(E) \simeq E$, we have an isomorphism 
$$
M_{n, 1}(E) \simeq E . 
$$

\begin{prop}\label{prop:srelliptic}
 Notation being as above, the map $\iota_*$ 
induces the isomorphism of schemes:
$$
\iota_*: M_{n, 1}(E) \simeq M_{n \cdot E}.
$$
Hence 
\begin{equation} \label{eq:isom}
M_{n \cdot  E} \simeq E. 
\end{equation}
and  $M_{n \cdot E}$ is a smooth irreducible component 
of $M_{n \cdot [E]}(X)$. 
\end{prop}

\vspace{0.5cm}
{\it Proof.} Since $\iota_*$ is an injective morphism, 
we only have to prove the surjectivity of $\iota_*$. 
Let $\E$ be a stable $\O_X$-coherent sheaf of pure dimension $1$ 
with $supp(\E) = E$ and 
$P(\E, m) =\chi(\E \otimes L^m) = n d m + 1$.  Then we have to prove 
that there exists a $\O_C$-coherent sheaf 
$\F \in M_{n, 1}(E) $ such that 
\begin{equation}\label{eq:push}
\E = \iota_*(\F).
\end{equation}
We prove the claim above (\ref{eq:push}) by induction with respect to $n$.  
For $n = 1$, let $\E$ be a $\O_X$ coherent sheaf in $M_{n \cdot E}$.  
Since $\dim H^0(X, \E) \geq  \chi(\E) = P(\E, 0) = 1$, we have a
 non-trivial section $s \in H^0(X, \E)$ which defines a non-trivial 
 homomorphism 
 $$
s: \O_X \longrightarrow \E.
 $$
Then setting $J = \ker s$, we have an inclusion of $\O_X$-coherent sheaves: 
$$
\varphi_s :\O_X/J \hookrightarrow \E.
$$
Let $I= I_E$ be the ideal sheaf of $E \subset X$.  
In the same argument as in the proof of Proposition \ref{prop:srrational}, one can shows that $J = I$  or for some $k \geq 1$ there exists 
a surjection 
$$
I^k/I^{k+1} = S^k( L \oplus L^{-1}) \twoheadrightarrow \O_X/J.
$$
In each case, one can obtain  an exact sequence of $\O_X$-sheaves 
\begin{equation} 
0 \rightarrow \iota_*(\cG) \rightarrow \E \rightarrow \E/\iota_*( \cG) \rightarrow 0 
\end{equation}
where $\cG$ is a line bundle on $E$ of degree $0$.  
The additivity of Hilbert polynomial  $ P(\E, m) = P(\cG, m) + P(\E/\iota_*(\cG), m)$ with $ P(\cG, m) = d m $ implies that 
$ P(\E/\iota_*(\cG), m) = 1 $.  
Hence there exists a closed point $x \in E$ such that 
$$
\E/\iota_*(\cG) \simeq \O_{X, x}/ m_x \simeq \C_x.
$$  
(Here $m_x$ is the maximal ideal sheaf of $x$.) 
Therefore one obtains the exact sequence of $\O_X$--sheaves 
\begin{equation} \label{eq:ext-1}
0 \rightarrow \iota_*(\cG) \rightarrow \E \rightarrow \C_x \rightarrow 0 
\end{equation}
Next we show that 
\begin{equation}\label{eq:ext-2}
\mbox{\rm Ext}^1_{\O_X}(\C_x, \iota_*(\cG)) \simeq \mbox{\rm Ext}^1_{\O_E}( \C_x, \cG).
\end{equation}
If one can show (\ref{eq:ext-2}), we can conclude that 
the extension $\E$  of $\O_X$-coherent sheaves  can be written as $\E = \iota_*(\F_1)$  where $\F_1$ is an  $\O_E$-coherent sheaf on $E$ with an extension of $\cG$ and $\C_x$.   
It is easy to see that 
$ \F_1 \in M_{1, 1}(E)$ and this proves our claim (\ref{eq:push}) for $n=1$.  

Next let us show the assertion \eqref{eq:ext-2}. 

We consider the following 
exact sequence followed from the local to global 
spectral sequence.
\begin{equation}
\begin{array}{l}
0 \rightarrow 
 H^1( \cHom_{\O_X}(\C_x, \iota_* (\cG)))  \rightarrow  
\mbox{\rm Ext}^1_{\O_X}(\C_x, \iota_*(\cG))   \rightarrow  H^0( \cExt^1_{\O_X}(\C_x, \iota_*(\cG))  \\ 
   \quad \quad  \rightarrow 
  H^2( \cHom_{\O_X}(\C_x, \iota_* (\cG))) 
  \end{array}
\end{equation}
One can easily see that 
$$
\cHom_{\O_X}(\C_x, \iota_*(\cG)) \simeq \cHom_{\O_E} ( \C_x, \cG) = 0
$$
$$
 \cExt^1_{\O_X}(\C_x, \iota_*(\cG))  \simeq \cExt^1_{\O_E}(\C_x, \cG) \simeq \C_x 
$$
Then this implies that 
$$
\mbox{\rm Ext}^1_{\O_X}(\C_x, \iota_*(\cG)) \simeq 
H^0 (E, \cExt^1_{\O_E}(\C_x,\cG)) \simeq  
\mbox{\rm Ext}^1_{\O_E}(\C_x, \cG )  \simeq H^0(E, \C_x)  \simeq \C. 
$$
which proves the claim (\ref{eq:ext-2}) and completes the proof for 
$n =1$.  

Next we assume that the assertion (\ref{eq:push}) is true for $\E'
\in M_{k \cdot E}$  $ 1 \leq k \leq n-1$. 
 Take $\E \in M_{n \cdot E}$ with $n \geq 2$. 
 Then by the same argument
 there exists an invertible sheaf $\cG$ of degree $0$ on $E$ and an 
 injective homomorphism $\iota_*(\cG)    \subset \E $.   Set $\E' = \E/\iota_*(\cG)$.  Since $P(\cG, m) = m d $ and 
 $ P(\E', m) = n d m + 1 $, we see that $ P(\E', m) = (n-1) d m + 1$.  Moreover again from the stability of $\E$ it is easy to see that $\E'$ is a stable $\O_X$-coherent  sheaf of pure dimension $1$  with support $E$ if $n \geq 2$.  
 Therefore 
 $\E' \in M_{(n-1) \cdot E}$.  By the assumption of induction, we 
 have $\cF' \in M_{n-1, 1}(E)$ such that $ \E' = \iota_*(\cF')$.

%%% 
Next we see that 
 \begin{equation}\label{eq:ext-3}
\mbox{\rm Ext}^1_{\O_X}(\iota_*(\cF'), \iota_*(\cG)) \simeq \mbox{\rm Ext}^1_{\O_E}( \cF', \cG).
\end{equation}
If \eqref{eq:ext-3} is true, the extension $\E$ of $\iota_*(\cF')$ and $\iota_*(\cG)$ can be written as in (\ref{eq:push}) and 
this completes the proof.

In order to prove (\ref{eq:ext-3}), we again use the following exact sequence 
\begin{equation}\label{eq:ext-4}
\begin{array}{l}
0 \rightarrow 
 H^1( \cHom_{\O_X}(\iota_*(\cF'), \iota_* (\cG)))  \rightarrow  
\mbox{\rm Ext}^1_{\O_X}(\iota_*(\cF'), \iota_*(\cG))   \rightarrow  H^0( \cExt^1_{\O_X}(\iota_*(\cF'), \iota_*(\cG))  \\ 
   \quad \quad  \rightarrow 
  H^2( \cHom_{\O_X}(\iota_*(\cF'), \iota_*(\cG))) . 
  \end{array}
\end{equation}
Again we easily see that 
\begin{equation}\label{eq:loc-hom}
\cHom_{\O_X}(\iota_*(\cF'), \iota_* (\cG)) \simeq \cHom_{\O_E} ( \cF', \cG) .
\end{equation}
From the exact sequence 
$$
 0 \lra I \lra \O_X \lra \O_E \lra 0,  
$$
and the isomorphism  $ \cHom_{\O_X}(\O_E, \O_E) \stackrel{\simeq}{\lra}  \cHom_{\cO_X}(\O_X, \O_E) $, we obtain the isomorphism  
$$
\cHom_{\O_X}(I, \O_E) \simeq  \cExt^1_{\O_X}(\O_E, \O_E). 
$$
Since $\cHom_{\O_X}(I, \O_E) \simeq \cHom_{\O_E}(I/I^2, \O_E)$, we finally obtain the isomorphism 
\begin{equation}
\cExt^1_{\O_X}(\O_E, \O_E) \simeq \cHom_{\O_E}(I/I^2, \O_E) \simeq N 
\end{equation}
where $N = N_{E/X}$ is the normal sheaf for $E \subset X$. 

Moreover we have the isomorphism 
\begin{equation}\label{eq:loc-ext}
 \cExt^1_{\O_X}(\iota_*(\cF'), \iota_*(\cG))  \simeq 
 \cExt^1_{\O_X}(\O_E, \O_E) \otimes \cHom_{\O_E}(\cF', \cG) 
  \simeq N \otimes \cHom_{\O_E}(\cF', \cG).   
\end{equation}
(Note that the first isomorphism follows from the locally freeness of 
$\cF'$ and $\cG$.) 
Since $\cF'$ is a stable $\cO_E$-sheaf with $\deg \cF' = 1 > 0$ and $ \deg \cG \otimes N = 0$, we can conclude that  
\begin{equation}\label{eq:ext-van}
H^0(E, N \otimes \cHom_{\O_E}(\cF', \cG)) = Hom(\cF', \cG \otimes N) = \{ 0 \}. 
\end{equation}

%%%%%%% PROOF OF THE ABOVE FACT %%%%%%%%
%If this is not the case, we have a non-trivial 
%$\O_E$-homomorphism 
%$\phi: \cF' \lra \cG \otimes N $.  
%Since $ N = L \oplus L^{-1}$ with $
%\deg L = 0$ and $\deg \cG = 0$, 
%we see that $ \deg \im (\phi)  \leq 0$.  
%Note that $\deg \cF' = \deg %\im \phi + \deg \ker \phi$.  
%Then $ker \phi \subsetneqq \cF'$ is a locally free $\O_E$-sheaf with %degree $ \geq \deg \cF' = 1$,  which contradicts the stability of 
%$\cF'$.

From (\ref{eq:loc-hom}), (\ref{eq:loc-ext}) and  (\ref{eq:ext-van}), 
the sequence (\ref{eq:ext-4}) gives the isomorphism 
\begin{equation}
\Ext^1_{\O_X}( \iota_*(\cF'), \iota_*(\cG)) \simeq H^1(E,  \cHom_{\O_E} ( \cF', \cG)) \simeq \Ext^1_{\O_E}(\cF', \cG)
\end{equation}
as required in (\ref{eq:ext-3}).  

We have proved that $i_*: M_{n,1}(E) \lra M_{n \cdot E}$ gives the 
isomorphism of set of $\C$-valued points.  
Since $M_{n,1}(E) \simeq E$ is smooth one dimensional scheme, 
if we prove the Zariski tangent space at each point of 
$M_{n \cdot E}$ 
is one dimensional, $i_*$ gives an isomorphism of schemes.  
For each element $\E \in M_{n \cdot E}$, we have a locally 
free $\O_E$-sheaf $\cF$ of rank $n$ of degree $1$ 
such that $ \E = \iota_* (\cF)$.   
The Zariski tangent space $T_{[\E]}$ at 
$\E$ in $ M_{n \cdot E}$ is given by 
$$
T_{[\E]} = \Ext^1_{\O_X} (\E, \E) = 
 \Ext^1_{\O_X}(\iota_*(\cF), \iota_*(\cF)). 
$$
We shall prove 
\begin{equation}\label{eq:ext-6}
\Ext^1_{\O_X}(\iota_*(\cF), \iota_*(\cF)) \simeq \Ext^1_{\O_E}(\cF, \cF) \simeq \C.
\end{equation}
(The last isomorphism follows from the fact that 
 $M_{n,1}(E)$ is smooth  and one dimensional.)
From the similar exact sequence as (\ref{eq:ext-4}), if we show 
\begin{equation}\label{eq:ext-van1}
H^0( X, \cExt^1_{\O_X}(\iota_*(\cF), \iota_*(\cF))) = 0, 
\end{equation}
we have the isomorphism 

\begin{equation}
\begin{array}{l}
\Ext^1_{\O_X}(\iota_*(\cF), \iota_*(\cF)) \simeq
 H^1(X, \cHom_{\O_X}(\iota_*(\cF), \iota_*(\cF))) 
 \simeq  \\
 H^1(E, \cHom_{\O_E} ( \cF, \cF)) 
\simeq \Ext^1_{\O_E}(\cF, \cF). 
\end{array}
\end{equation}
as required.  For (\ref{eq:ext-van1}), we again 
recall  the isomorphism (cf. (\ref{eq:loc-ext}))
$$
\cExt^1_{\O_X}(\iota_*(\cF), \iota_*(\cF)) \simeq 
\cExt^1_{\O_X}(\O_E, \O_E) \otimes \cHom_{\O_E}(\cF, \cF) \simeq N \otimes  \cHom_{\O_E}(\cF, \cF), 
$$
which induces the isomorphism 
$$
H^0( X, \cExt^1_{\O_X}(\iota_*(\cF), \iota_*(\cF))) \simeq H^0(E, N \otimes  \cHom_{\O_E}(\cF, \cF) ). 
$$
Since $N = L \oplus L^{-1}$,  the non-vanishing of  (\ref{eq:ext-van1}) implies the existence of  a non-trivial isomorphism  of $\O_E$-sheaf 
$$
\phi^+: \cF \lra \cF \otimes L \quad  \mbox{or} \quad  \phi^-: \cF \lra \cF \otimes L^{-1}.
$$
Since $\deg L = 0 $, the stability of $\cF$ again implies that 
$\phi^{\pm}$ is an isomorphism if it is non-trivial.  Taking the 
highest power  of $\phi^{\pm}$, we obtain the isomorphism
$$
\wedge^{n} \phi^{\pm}: \wedge^n \cF \stackrel{\simeq}{\lra}  
(\wedge^n \cF) \otimes L^{\pm n}.  
$$ 
This  implies that $L^{\pm n} \simeq \O_E$, which contradicts to 
the fact that $L$ is a non-torsion element of $\Pic(E)$.  
Thus we have the vanishing (\ref{eq:ext-van1}),  and hence (\ref{eq:ext-6}).  

\qed

By Proposition \ref{prop:srrational} and \ref{prop:srelliptic}, we can calculate our local BPS invariants as follows.  Note that 
our local BPS invariants  coincide with the conjectural local BPS 
invariants predicted by Conjecture  \ref{conj:gvc} and the calculations of Gromov--Witten invariants (cf. Theorem \ref{thm:srr} 
and \ref{thm:sre}).

\begin{prop}\label{prop:sr-bps}  
Let $C \simeq \BP^1  \subset X$ be a rigid rational curve in a 
smooth Calabi-Yau 3-fold.  Then our local BPS invariants for $C$ are given as follows: 
\begin{equation}\label{eq:sr-bps}
n_g( n \cdot \BP^1) = 
\left\{\begin{array}{l}
1 \quad \mbox{for  $g=0$ and  $n  = 1$} \\
0 \quad \mbox{otherwise. }
\end{array} \right. 
\end{equation}
\end{prop}

\begin{prop}\label{prop:se-bps}
  Let $E \simeq \BP^1  \subset X$ be a super-rigid elliptic curve in a smooth Calabi-Yau 3-fold.  Then our local BPS invariants for $C$ are given as follows: 
\begin{equation}\label{eq:sr-bpse}
n_g( n \cdot E) = 
\left\{\begin{array}{l}
1 \quad \mbox{for  $g=1$ and all $n \geq 1$} \\
0 \quad \mbox{otherwise. }
\end{array} \right. 
\end{equation}
\end{prop}

\vspace{0.5cm}
\noindent
{\it Proof of Proposition \ref{prop:sr-bps} and \ref{prop:se-bps}. }

\vspace{0.5cm}
For the case of $C = \BP^1$, $M_{ 1 \cdot \BP^1} \simeq \Spec \C$ for $n =1$ and $M_{n \cdot \BP^1}$ is empty if $ n > 1$.  Then 
the $sl_2 \times sl_2$ decomposition of 
$$
IH^*( M_{ 1 \cdot \BP^1}) = (0)_L \otimes (0)_R \simeq (I_1)^{0} \otimes (0) . 
$$
Hence $R_0 ( 1 \cdot \BP^1) = (0)_R $  and $R_h (1 \cdot \BP^1) = \emptyset $ for $ h > 0$.  Therefore  we have 
$$
n_0( \BP^1) =  Tr_{R_0 ( 1 \cdot \BP^1)}(- 1)^{2 H_R} =  (-1)^0 = 1. 
$$
Moreover $ R_h(n \cdot \BP^1) = \emptyset $ unless $ n = 1$ and $h = 0$.  Hence we have the assertion for $C = \BP^1$.  

Next let us consider the case of a super-rigid 
elliptic curve $E \subset X$.  
In this case, for all integer $n \geq 1$, 
  $M_{n \cdot E}(X)$ is always isomorphic to $E$. 
Moreover the image of the  support map 
$$
\pi: M_{n \cdot E}(X) \lra Chow_{[n \cdot E]}
$$
is just one point $\{ n \cdot E \}$.  
Hence we can identify it with the 
structural morphism  $ \pi: M_{n \cdot E} \lra  \Spec \C$.  
Moreover from Proposition \ref{prop:srelliptic}, we have 
$$
M_{n \cdot E} \simeq E.  
$$ 
Hence the $sl_2 \times sl_2$ decomposition for the morphism $\pi$ 
is given by 
$$
IH^*(M_{n \cdot E}) \simeq H^*(E, \C) \simeq I_1  \otimes (0)_R .
$$
Therefore we have $R_1( n \cdot E) \simeq (0)_R$ and $R_h ( n \cdot E)$ is empty for $ h \not= 1$.  Therefore 
$$
n_1( n \cdot E) = 1 
$$
and 
$$
n_h (n \cdot E) =  0  \quad \mbox{otherwise}. 
$$

\qed

%%%%%%%%%%%%%%%%%%%%%%%%%%%%%%%%%%%%%%%%%%%%%%%%%%%%%%%%%%%%%%%%%%%%%%%%%%%%%%

\vspace{0.5cm} 
\subsection{Rational Elliptic Surfaces in Calabi--Yau $3$-folds}
~~~

\vspace{0.5cm}
Next, we shall calculate our local   BPS invariants for 
special homology class of a rational elliptic surface $S$  in a smooth Calabi--Yau 3-fold $X$ based on our mathematical Definition \ref{def:bps}.  We should remark that in \cite{hst:1} by a physical argument (holomorphic anomaly equation) and explicit 
calculations of some part of Gromov--Witten invariants we calculated 
a part of the left hand side of \eqref{eq:gv}.  Then using Jacobi triple product formula, we found that the  
corresponding right hand side of the conjectural 
formula \eqref{eq:gv} 
can be  obtained from the $sl_2 \times sl_2$ spin 
version of G\"{o}ttsche's formula for the Poincar\'e polynomial 
of the Hilbert schemes $S^{[g]}$ of 
$g$-points on a rational elliptic surface $S$.  

In what follows, we shall  prove that the above physical arguments 
can be justified in our mathematical definition of BPS invariants in 
\eqref{def:bps} using the moduli of sheaves.

Let $f:S \lra \BP^1$ be a rational elliptic surface with a section 
$\sigma:\BP^1 \lra S$,  and  $C = \sigma(\BP^1)$ the image of the section.  (Note that $C$ is a smooth rational curve with $C^2 = -1$ and hence  $C$ is rigid in $S$. )  
We denote by $ F $ the divisor  class of a fiber of $f$.     
We  consider the case when $S$ lies in a smooth Calabi-Yau 3-fold $X$ as a divisor. (For example, see \cite{hss}, \cite{saito}).   
For a non-negative integer $g$, consider the homology class 
\begin{equation}\label{eq:sec}
\beta_g =[ C + g F]  \in  H_2 (X, \Z). 
\end{equation}

In order to consider the  BPS invariant 
$n_h(\beta_g) = n_h([C + g F])$, 
we have to take  account of 
all curves which are homological 
equivalent to $\beta_g$.  
Instead, we will consider the local 
BPS invariants $n_h( C + g F)$ which 
are defined by the moduli of sheaves 
whose support are curves in the rational surface $S$.  
Since $C \subset S \subset X$ is a rigid rational 
curve both in $X$, this makes sense.  

More precisely, let us consider the 
linear system $| C + g F|$ which consists of effective divisors 
linearly equivalent to the divisor $C + g F$ of $S$.  
Then it is easy to see that $| C + g F| \simeq C + |g F| \simeq \BP^g$.    

Let us fix an ample line bundle $ L $ on $X$. 
For simplicity, we assume that $L$ 
can be chosen as\footnote{We do not know that this assumption is always true, but we know examples of $S \subset X$  which satisfies the condition \eqref{eq:ample}. See \cite{hss}. }
\begin{equation}\label{eq:ample}
L_{|S} = C + k F, \quad k >> 0.  
\end{equation}
and set $d = L \cdot (C + g F)$.   

Then 
the moduli space which we should  consider is given by 
\begin{equation}
M_{C + g F} := \left\{ \begin{array}{c|l}  
 & \mbox{a semistable sheaf on $X$, pure of  } \\
\ \E \  & \mbox{dimension $1$ with} \  P(\E, m) = d m + 1 \\
  & and   \  s(\E) \in |C + g F|  \end{array}     \right\}/ \mbox{isom.} 
\end{equation}
Then we obtain  the natural support  map
\begin{equation}
\begin{array}{ccc}
\pi:M_{C + g F} & \lra &   |C + g F | \\
      \E  & \mapsto &  \pi(\E) := s(\E).   
      \end{array} 
\end{equation}
The local BPS invariant $n_h(C + g F)$ is defined by using the relative 
Lefschetz decomposition 
$$
IH^*(M_{C+ g F}, \C)  = \bigoplus_{h \geq 0 } I_h \otimes R_h(C + g F)
$$
as
\begin{equation}
n_h(C + g F) := Tr_{R_h(C + g F)} (- 1)^{2H_R}.
\end{equation}

The following is our  main theorem (cf. \cite{hst:1}).  
\begin{thm}\label{thm:main}
\begin{multline*}
\sum_{g \geq 0, h \geq 0}n_h(C+gF)\left(2\sin\frac{\lambda}{2}\right)^{2h-2} q^g \\
=\frac{1}{(e^{-\sqrt{-1}\lambda/2}-e^{\sqrt{-1}\lambda/2})^2}
\prod_{n\ge 1}\frac{1}{(1-e^{-\sqrt{-1}\lambda}q^n)^2(1-e^{\sqrt{-1}\lambda}q^n)^2
(1-q^n)^8}.
\end{multline*}
\end{thm}
Fix an ample divisor $L$ as in \eqref{eq:ample}. 
Let $M^{ss}(r,c_1,\chi)$ be the moduli space of semistable sheaves $\E$ on $S$  with rank $r(\E)=r$, $c_1(\E)=c_1$ and $\chi(\E)=\chi$ with respect to $L$.  Moreover we assume that all closed 
fibers  of rational elliptic surface $f:S \lra \BP^1$ are integral.   
The following proposition is due to Kota Yoshioka (cf. \ \cite{yo:1}).

\begin{prop} \label{prop:yoshioka}
Under the assumption as above, we have the following. 
\begin{enumerate} 
\item We have an isomorphism 
 $M^{ss}(0, C + g F, 1) \simeq M_{C + g F} $. 
\item Fourier--Mukai transform induces an isomorphism
$$
M^{ss}(0, C+gF ,1) \simeq M^{ss}(1, 0,1-g) \simeq S^{[g]},
$$
where $S^{[g]}:=Hilb^g(S)$ is the Hilbert scheme of $g$ points on $S$.
Hence we have an isomorphism 
$$
 M_{C + g F} \simeq S^{[g]}.  
$$
\item Moreover the support map $\pi: M_{C + g F} \lra | C + g F| \simeq \BP^g $ can be identified with the natural map 
\begin{equation}\label{eq:support}
\pi:S^{[g]} \lra Sym^g (\BP^1) \simeq \BP^g 
\end{equation}
given by the composite map of  Hilbert--Chow morphism 
$ S^{[g]} \lra Sym^g (S)$ and the natural map $ Sym^g (f): Sym^g (S) \lra Sym^g(\BP^1) $.
\end{enumerate}
\end{prop}

From Proposition \ref{prop:yoshioka}, in order to show 
Theorem \ref{thm:main},  we only have to determine the 
relative Lefschetz decomposition of $IH^*(S^{[g]}, \C) \simeq H^*(S^{[g]}, \C) $ with respect to the natural morphism 
\begin{equation}\label{eq:mor}
\pi:S^{[g]} \lra \BP^{g}.
\end{equation}

We recall here  well-known G\"ottsche's formula for 
the Poincar\'e polynomial of the Hilbert schemes $Y^{[g]} = Hilb^g(Y)$ for a smooth projective surface $Y$.
Let $b_i(Z)$ be the $i$-th Betti number of a smooth projective variety $Z$.  Let us define the shifted Poincar\'e polynomial  by  
$$
P_t(Z) =t^{- \dim Z} \cdot \left( \sum_{ k \geq 0}  b_k(Z) t^k   
\right),  
$$ 
where $\dim Z$ is the complex dimension of $Z$.  
Note that $P_{t^{-1}}(Z) = P_t(Z)$.  
The following formula is proved by G\"ottsche \cite{got}.  
For a proof see \cite{got}, \cite{got-soe} or \cite{nak}.

\begin{thm}  For a smooth projective surface $Y$, 
the generating function of shifted 
Poincar\'e polynomials of the Hilbert schemes $Y^{[g]} $ of 
$g$-points is given by the formula
\begin{equation}\label{eq:got}
\sum_{g\ge 0}P_t(Y^{[g]})q^g=\prod_{n\ge 1}
\frac{(1 + t^{-1} q^n)^{b_1(Y)}(1 + t \ q^n)^{b_3(Y)}}{(1-t^{-2}q^n)^{b_0(Y)}(1-t^2q^n)^{b_4(Y)}(1-q^n)^{b_2(Y)}}.  
\end{equation}
\end{thm}

For a rational elliptic surface $ f : S \lra \BP^1$, 
it is easy to see that $b_1(S) = b_3(S) = 0$ and  $ b_2(S) = 10$.  Therefore the formula 
\eqref{eq:got} is  reduced to 
$$
\sum_{g\ge 0}P_t(S^{[g]}) q^g= \prod_{ n \geq 1} \frac{1}{
(1-t^{-2}q^n)(1-t^2q^n)(1-q^n)^{10}}.  
$$
Recall that an ample line bundle $L'$ on $S^{[g]}$ defines 
the usual Lefschetz $sl_2$ action on the cohomology ring 
$H^*(S^{[g]}, \C)$.  Denote by $H$ the usual weight operator for 
the $sl_2$ action.  Then the Poincar\'e polynomial 
can be expressed as 
the character of  the  representation:
$$
P_t(S^{[g]})=Tr_{H^*(S^{[g]},\C)}(-1)^{2H} t^{2 H }. 
$$

The relative Lefschetz action on  $H^{*}(S^{[g]}, \C)$ 
with respect to the morphism $\pi: S^{[g]} \lra \BP^g$ (cf. \eqref{eq:mor})  determines the 
representation of $(sl_2)_L\times (sl_2)_R$ and let us consider its 
character
$$
P_{t_L,t_R}(S^{[g]})=Tr_{H^*(S^{[g]},\C)}(-1)^{2H_L+2H_R}t_L^{2H_L}t_R^{2H_R}.
$$

We have the following generalization of G\"ottsche's formula.  

\begin{thm} \label{thm:gen-got}
Let $f:S \lra \BP^1$ be a  
rational elliptic surface. Assume that all closed fibers are 
integral.  Then we have 
\begin{multline}\label{eq:gen-got}
\sum_{g\ge 0}P_{t_L,t_R}(S^{[g]})q^g=\prod_{n\ge 1}
\frac{1}{(1-(t_Lt_R)^{-1}q^n)(1-(t_Lt_R)q^n)}\\
\times\prod_{n\ge 1}\frac{1}{(1-(t_Lt_R^{-1})q^n)(1-(t_L^{-1}t_R)q^n)
(1-q^n)^{8}}.
\end{multline}
\end{thm}

For a proof of Theorem \ref{thm:gen-got}, we need the following
\begin{prop} \label{prop:sym}
Let $f:S \lra \BP^1$ be a  
rational elliptic surface. Assume that all closed fibers are 
integral\footnote{This condition is just for simplicity.  Even if we do not assume that all closed fibers are integral, we can obtain the  statement. See remark below.} . 
\begin{enumerate}

\item Set $E^{r,s} := {}^p H^r( \BP^1,  {}^p R^{s} f_* \C_S [2] )$. Then 
we have isomorphisms:
\begin{center}
\begin{tabular}{|ccc} 
 $E^{1, -1}$ & $E^{1, 0}$ & $ E^{1, 1}$  \\
$E^{0, -1}$ & $E^{0, 0}$ &  $E^{0, 1}$  \\
$E^{-1, -1}$ & $E^{-1, 0}$ &  $E^{-1, 1}$ \\ \hline 
\end{tabular} 
\quad 
$=$ \quad \begin{tabular}{|ccc}
    $\C$ & $0$ & $\C$  \\
    $0$ & $\C^8 $ & $0$ \\
    $\C$ & $0$ & $\C$ \\ \hline 
 \end{tabular}. 
\end{center}
Moreover the $sl_2 \times sl_2$ decomposition of 
cohomology ring $H^*(S, \C[2])$ can be obtained as 
$$
H^*(S, \C[2]) \simeq (\frac{1}{2})_L \otimes (\frac{1}{2})_R +  
8 \cdot (0)_L \otimes (0)_R.  
$$ 
Therefore we have  
$$
 P_{t_L, t_R}(S) = (t_L t_R)^{-1}  + t_L (t_R)^{-1} 
+ (t_L)^{-1} t_R + t_L t_R +  8 .
$$

\item Let $Sym^{n}(f): Sym^{n}(S) \lra Sym^{n} (\BP^1) \simeq \BP^n$ be the natural morphism between $n$-th symmetric products 
of $S$ and $\BP^1$.   Let $P_{t_L, t_R}(Sym^n(S)) $ denote the character of $sl_2 \times sl_2$ representation of 
$IH^{*}(Sym^n(S), \C[2n]) \simeq H^*(Sym^n(S), \C[2n])$\footnote{Since $Sym^n(S)$ has only 
quotient singularities, the intersection cohomology groups are isomorphic to the ordianry 
cohomology groups}.
Then we have the following formula:
\begin{multline}\label{eq:formula}
\sum_{n \geq 0} P_{t_L,t_R}(Sym^n(S)) q^n =
\frac{1}{(1-(t_Lt_R)^{-1}q)(1-(t_Lt_R)q)}\\
\times \frac{1}{(1-(t_Lt_R^{-1})q)(1-(t_L^{-1}t_R)q)
(1-q)^{8}}.
\end{multline}
\end{enumerate}

\end{prop}

\vspace{0.5cm}
\noindent
{\it Proof of Proposition \ref{prop:sym}. } 
Since all fibers of $f$ are integral, we have 
\begin{equation}\label{eq:high-direct}
^p R^{-1} f_* \C_S [2] \simeq \C_{\BP^1}[1], \quad 
^p R^{1} f_* \C_S [2] \simeq \C_{\BP^1}[1].
\end{equation}
Hence,  we see that $E^{r, s} \simeq \C$ for $ (r,s) = (\pm 1, \pm 1)$.  Moreover 
since the Leray spectral sequence degenerates and $ \dim E^{0,0} + \dim E^{-1, 1} + \dim E^{1, -1} =  
\dim  \H^0(S, \C_S[2]) = \dim H^2(S, \C) = 10$, we have $\dim E^{0,0} = 8$.  (Note that $E^{-1, 1} \simeq 
H^0(\BP^1, R^2 f_* \C) \simeq \C $ is the space of the class of a fiber  of $f$ and  
$E^{1, -1} = H^2(\BP^1, R^0 f_* \C) \simeq \C $ is the space of the class of the section).  Therefore, 
the relative Lefschetz action and hence $sl_2 \times sl_2$ decomposition  
of $H^*(S, \C[2])$ 
with respect to $f:S \lra \BP^1$ are  
determined by this  Leray spectral sequence.  
For the $sl_2 \times sl_2$-decomposition of the intersection cohomology group 
$IH^*(Sym^n(S), \C[2n]) \simeq H^*(Sym^n(S), \C[2n])$ of 
 the symmetric power $Sym^n(S)$, we recall that 
 $ IH^*(Sym^n(S), \C) \simeq H^*(Sym^n(S), \C) $ is ${\mathfrak S}_n$-invariant part of $H^*(S^n, \C) = (H^*(S, \C))^{\otimes n}$. 
  Moreover the $sl_2 \times sl_2$-decomposition of $(H^*(S, \C))^{\otimes n}$ 
 is nothing but the one induced by $H^*(S,\C)$.  Hence we have the 
 formula  \eqref{eq:formula}. 
\qed

\begin{rem}
In the proof above, we do not have to 
assume that all fibers of $f:S \lra \BP^1$ are integral.   
First we allways have ${}^p R^{-1} f_* \C_{S}[2] \simeq R^0 f_* \C_S \simeq \C_{\BP^1}$ (up to a shift) by the connectivity of 
fibers.  Then by 
the relative Lefschetz theorem we have the isomorphism 
$ {}^p R^{-1} f_* \C_{S}[2] \simeq {}^p R^{1} f_* \C_{S}[2]$.  
This shows that the isomorphism \eqref{eq:high-direct} is still true, and so are the all assertions of Proposition \ref{prop:sym}.  
Note that,  if  some fibers of $f$ are reducible, the  
 usual higher direct image sheaf $ R^{2} f_* \C_S $ becomes a direct 
 sum of $\C_{\BP^1}$ and skyscraper sheaves supported on 
 points corresponding to reducible fibers.  Hence for usual higher 
 direct image sheaves, the relative hard Lefschetz theorem
 $ R^0 f_* \C_S \simeq R^{2} f_* \C_S$ does not 
 holds if some fibers of $f$ are reducible, contray to the case of 
 perverse higher direct image sheaves.   

\end{rem}

\vspace{0.5cm}
\noindent
{\it Proof of Theorem \ref{thm:gen-got}.}

In [6.2, \cite{nak}], Nakajima gives a proof of G\"ottsche's formula 
using the perverse sheaf and the fact that the Hilbert--Chow 
morphism $S^{[n]} \lra Sym^n(S)$ is semismall.  (See also \cite{got-soe}.)  

Consider the Hilbert--Chow morphism $S^{[n]}\to {\rm Sym}^n(S) $ and 
stratification of 
$$
{\rm Sym}^n(S)= \bigcup_\nu {\rm Sym}^n_\nu S,
$$
$$
{\rm Sym}^n_\nu S:=\{\sum_{i=1}^k\nu_i[x_i]\in{\rm Sym}^gS~|~x_i\ne x_j~
{\rm for~all}~i\},~~~\nu_1\ge\nu_2\ge\dots\ge\nu_k,
$$
defined by the partitions $\nu$ of $n $. 
For a partition $\nu$, define $ \alpha_i:= \sharp \{l|\nu_l=i\}$. 
Then we can define 
$$
{\rm Sym}^{\nu}(S) := {\rm Sym}^{\alpha_1}S\times\dots\times{\rm Sym}^{\alpha_n} S.  
$$
Let $2l(\nu)$ be the complex dimension of ${\rm Sym}^{\nu}(S)$.  
Then one can see (cf. [(6.13), \cite{nak}]):
\begin{equation}\label{eq:decomp}
H^{i + 2n}( S^{[n]}, \C) = \oplus_{\nu} H^{i + 2l( \nu)}({\rm Sym}^\nu(S), \C).
\end{equation}
Since   $sl_2 \times sl_2$-decompositions of both side of 
\eqref{eq:decomp} are compatible with each other, 
 by the same argument as in [6.2, \cite{nak}] together with the formula \eqref{eq:formula}, we can show  the formula  \eqref{eq:gen-got}. 

\qed

\noindent
{\it Proof of Theorem \ref{thm:main}. }

Now, we shall prove Theorem \ref{thm:main}.  

Let 
$$
H^*(S^{[g]},\C) := \oplus_{ h \geq 0} I_h \otimes R_h(C + g F) 
$$
be the $sl_2 \times sl_2$-decomposition of  $H^*(S^{[g]},\C)$. 
In the left action, the character of $I_1= (\frac{1}{2})_L + 2(0)_L$ is given by $(- t_L - t_L^{-1} +2)$ and hence the character of 
$I_h = (I_1)^{\otimes h}$ is given by $(- t_L - t_L^{-1} +2)^h$.  

Since the character of this decomposition is given by $P_{t_L t_R}(S^{[g]})$  in \eqref{eq:gen-got}, 
\begin{equation}
\begin{array}{ccc}
P_{t_L, t_R}(S^{[g]})_{|t_R =1} &  = & \sum_{h \geq 0} Tr_{R_h(C + g F)} (-1)^{2H_R} ( -t_L - t_L^{-1} + 2)^{h} \\
 &  = &    \sum_{h \geq 0} n_h(C + g F) ( -t_L - t_L^{-1} + 2)^{h}
\end{array}
\end{equation}
Then setting  $t_L = e^{\sqrt{-1} \lambda}$, we have 
$$
(- t_L - t_L^{-1} + 2) = -e^{\sqrt{-1} \lambda}- e^{-\sqrt{-1} \lambda} + 2 = \left( 2 \sin (\frac{\lambda}{2})\right)^2. 
$$
From Theorem \ref{thm:gen-got}, we obtain 
\begin{align*}
&\sum_{g\ge 0, h \ge 0}n_h(C+gF)\left(2\sin\frac{\lambda}{2}\right)^{2h-2}q^g\\
=&\frac{1}{(e^{-\sqrt{-1}\lambda/2}-e^{\sqrt{-1}\lambda/2})^2}
\sum_{g\ge 0}P_{t_L=e^{\sqrt{-1}\lambda},t_R=1}(S^{[g]})q^g\\
=&\frac{1}{(e^{-\sqrt{-1}\lambda/2}-e^{\sqrt{-1}\lambda/2})^2}
\prod_{n\ge 1}\frac{1}{(1-e^{-\sqrt{-1}\lambda}q^n)^2(1-e^{\sqrt{-1}\lambda}q^n)^2
(1-q^n)^8}.
\end{align*}
\qed

\vspace{0.2cm}
\begin{rem}
It is clear from the proof that Theorem \ref{thm:main}  also holds for other elliptic surfaces 
in a Calabi--Yau manifold. 
In particular, if we consider an elliptic K3 surface, we have the same results 
as that of Kawai--Yoshioka \cite{ky:1}. 
They considered the Abel--Jacobi map and counted the number of BPS states from 
D0-D2 system. 
They defined the moduli space of D0-D2 system as relative Hilbert schemes of 
$d$-points ${\mathcal C}_h^{[d]}\to |C_h|\simeq \P^h$, $d\ge 0$ where $C_h\subset K3$ is 
a curve of genus $h$ and ${\mathcal C}_h^{[1]}\simeq {\mathcal C}_h$ is the universal 
family over $|C_h|$.
$$
\sum_{h,d}\chi({\mathcal C}_h^{[d]})q^hy^{d+1-h}=\frac{1}{(y^{1/2}-y^{-1/2})^2}
\prod_{n\ge 1}\frac{1}{(1-yq^n)^2(1-y^{-1}q^n)^2(1-q^n)^{20}}.
$$
On the other hand, we use the relative Lefschetz action on the relative Jacobian 
and counted the spin contents of BPS states from M2 brane. 
The coincidence of these results is very natural since the original physical 
theory is equivalent.
\end{rem}
If we allow some physical arguments (holomorphic anomaly equation), 
we have the nontrivial evidence of Gopakumar--Vafa conjecture.
Let us write the generating functions of 
Gromov--Witten invariants as
$$
Z_{g;n}(q):=\sum_{d}N_{g,d;n}q^d,~~~
N_{g,d;n}:=\sum_{(\beta,\sigma)=d,(\beta,F)=n}N_g(\beta),~~~n\ge 1,
$$
where $N_g(\beta) \in \Q$ are genus $g$ 
Gromov--Witten invariants for $\beta\in H_2(S,\Z)$ defined by 
$$
N_g(\beta):=\int_{\left[\overline{\M}_{g,0}(S,\beta)\right]^{virt}}
c_{{\rm top}}(R^1\pi_*\mu^*N_{S/X}).
$$
Explicitly, $Z_{0;1}(q)$ is given by (\cite{hss})
\begin{equation}\label{eq:z1}
Z_{0;1}(q)=E_4(q)\prod_{k\ge 1}\frac{1}{(1-q^k)^{12}}. 
\end{equation}

For $Z_{g;n}(q)$, we have the following conjecture in mathematics 
suggested  by some arguments in physics (cf. \cite{hst:1}).    
\begin{conj}$(${\bf Holomorphic anomaly equation \cite{hst:1}}$)$\\

\begin{enumerate}
\item 
 $Z_{g;n}(q)$ has the following expression
$$
Z_{g;n}(q)=\frac{P_{2g+6n-2}(E_2(q),E_4(q),E_6(q))}
{\prod_{k\ge 1}(1-q^k)^{12n}},
$$
where $P_{2g+6n-2}(E_2(q),E_4(q),E_6(q))$ is a homogeneous 
polynomial of weight $2g+6n-2$ and $E_*(q)$ are Eisenstein 
series of weight $*$.

\item 
$P_{2g+6n-2}(E_2,E_4,E_6)$ satisfies the following equation:
\begin{multline*}
\frac{\partial P_{2g+6n-2}}{\partial E_2}=\frac{1}{24}
\sum_{g=g'+g''}\sum_{s=1}^{n-1}s(n-s)P_{2g'+6s-2}P_{2g''+6(n-s)-2}\\
+\frac{n(n+1)}{24}P_{2(g-1)+6n-2}.
\end{multline*}

\end{enumerate}
\end{conj}

\vspace{0.5cm}
We can solve the holomorphic 
anomaly equation easily (cf. \cite{hst:1})  and, if $n=1$,  
the generating function of $Z_{g;1}(q)$ may be summarized  to  
$$
\sum_{g\ge 0}Z_{g;1}(q)\lambda^{2g}=Z_{0;1}(q)\exp
\left(2\sum_{k\ge 1}\frac{\zeta(2k)}{k}E_{2k}(q)
\left(\frac{\lambda}{2\pi}\right)^{2k}\right).
$$
By the famous Jacobi's triple product formula, we have
\begin{multline*}
\lambda^{-2}\exp
\left(2\sum_{k\ge 1}\frac{\zeta(2k)}{k}E_{2k}(q)
\left(\frac{\lambda}{2\pi}\right)^{2k}\right)\\
=\frac{1}{(e^{-\sqrt{-1}\lambda/2}-e^{\sqrt{-1}\lambda/2})^2}\prod_{n\ge 1}
\frac{(1-q^n)^4}{(1-e^{\sqrt{-1}\lambda}q^n)^2(1-e^{-\sqrt{-1}\lambda}q^n)^2}.
\end{multline*}
Multiplying $Z_{0;1}(q)$ both sides, we can easily verify the Gopakumar--Vafa 
conjecture, which was given in our previous paper \cite{hst:1}.

Finally we remark that our horomorphic anomaly equation suffices to 
determine $Z_{g;n}(q)$ recursively for all $g$ and $n$, see [HST] for 
details. Here we present the first few solutions for $n=1,2$;
\begin{eqnarray*}
&&
 Z_{1,1}(q)=\frac{E_2(q) E_4(q)}{\prod_{n\geq1}(1-q^n)^{12}}  \;,\;\;
 Z_{2,1}(q)=\frac{E_4(q) ( 5 E_2(q)^2 + E_4(q) )}
            {1440 \prod_{n\geq1}(1-q^n)^{12}}   \\
&&
Z_{3,1}(q)=\frac{E_4(q)( 35 E_2(q)^3 + 21 E_2(q) E_4(q) + 4 E_6(q) )}
{362880 \prod_{n\geq1}(1-q^n)^{12}}  \;\;, 
\end{eqnarray*}

\begin{eqnarray*}
& Z_{0,2}(q)=&\frac{  E_2(q) E_4(q)^2 + 2 E_4(q) E_6(q)  }
                   { \prod_{n\geq1}(1-q^n)^{24} }  \;, \\  
& Z_{1,2}(q)=&\frac{  10 E_2(q)^2 E_4(q)^2 + 9 E_4(q)^3 + 
                      24 E_2(q) E_4(q) E_6(q) + 5 E_6(q)^2  }
               { 1152  \prod_{n\geq1}(1-q^n)^{24} }
                \;\;, \\
& Z_{2,2}(q) =& \big(   
                190 E_2(q)^3 E_4(q)^2 + 417 E_2(q) E_4(q)^3 + 
                540 E_2(q)^2 E_4(q) E_6(q) + \\   
&             &  356 E_ 4(q)^2 E_6(q) + 
                        225 E_2(q) E_6(q)^2  \big) 
                \frac{1} {207360 \prod_{n\geq1}(1-q^n)^{24} } 
             \;, \\
& Z_{3,2}(q)=& \big( 2275 E_2(q)^4 E_4(q)^2 + 8925 E_2(q)^2 E_4(q)^3 + 
                      3540 E_4(q)^4 +  \\
&              &7560 E_2(q)^3 E_4(q) E_6(q)  
                       + 14984 E_2(q) E_4(q)^2 E_6(q) +  \\
&              &4725 E_2(q)^2 E_6(q)^2 + 4071 E_4(q) E_6(q)^2  \big)
               \frac{1}{ 34836480 \prod_{n\geq1}(1-q^n)^{24} } \;\;.
\end{eqnarray*}

\vspace{0.5cm}
\begin{center}
{\bf Acknowledgements}
\end{center}

\vspace{0.2cm}
We would like to thank Kota Yoshioka for useful discussions about  
the moduli theory  of stable sheaves which are inevitable for finishing  section 3 and section 4.  We also thank 
R. Donagi, J. Bryan for useful discussions about BPS state 
countings in RIMS project ``Geometry related to string theory", Kyoto,  2000.  

%%%%%%%%%%%%%%%%%%%%%%%%%%%%%%%%%%%%%%%%%%%%%%%%%%%%%%%%%%%%%%%%%%%%%%%%%%%%%%
%%%%%%%%%%%%%%%%%%%%%%%%%%%%%%%%%%%%%%%%%%%%%%%%%%%%%%%%%%%%%%%%%%%%%%%%%%%%%%

%
%%%%%%%%%%%%%%%%%%%%%%%%%%%%%%%%%%%%%%%%%%%%%%%%%%%%%%%%%%%%%%%%%%%%%%%%%%%%%%
%
\end{document}